\documentclass[12pt,reqno]{amsart}
\usepackage{amsmath}
\usepackage{amssymb}
\usepackage{amsfonts}
\usepackage{mathrsfs}
\usepackage{graphicx}
\usepackage{epsfig}
\usepackage{enumitem}
\usepackage[T1]{fontenc}
\numberwithin{equation}{section}       
\usepackage{longtable}
\usepackage{threeparttable}

\usepackage{booktabs}
\usepackage[font=normalfont,tableposition=top,figureposition=bottom,skip=0pt]{caption}
\usepackage[utf8]{inputenc} 
\usepackage[T1]{fontenc} 
\usepackage{lmodern}
\usepackage{xcolor}

\theoremstyle{plain}
\newtheorem{theorem}{Theorem}[section]

\newtheorem{prop}{Proposition}[section]

\newtheorem{lemma}[prop]{Lemma}

\usepackage{hyperref}

\theoremstyle{definition}
\newtheorem{definition}[prop]{Definition}

\theoremstyle{remark}

\newtheoremstyle{citing}
  {3pt}
  {3pt}
  {\itshape}
  {}
  {\bfseries}
  {.}
  {.5em}
  {\thmnote{#3}}

\theoremstyle{citing}

\DeclareMathAlphabet{\mathpzc}{OT1}{pzc}{m}{it} 

\newcommand{\hyp}{\textrm{hyp}}

\newcommand{\N}{\mathbb{N}}

\newcommand{\Q}{\mathbb{Q}}
\newcommand{\R}{\mathbb{R}}
\newcommand{\T}{\mathbb{T}}

\newcommand{\Z}{\mathbb{Z}}

\newcommand{\teta}{\widetilde{\teta}}

\newcommand{\eps}{\varepsilon}

\newcommand{\dist}{d}

\DeclareMathOperator{\diam}{diam}

\DeclareMathOperator{\supp}{supp} 

\def\marginparQ#1{\marginpar{}}
\begin{document}
\title[]{Multifractal Analysis of generalized Thue-Morse trigonometric polynomials}

\author{Aihua FAN}
\address{
	LAMFA, UMR 7352 CNRS, University
	of Picardie, 33 rue Saint Leu, 80039 Amiens, France}
\email{ai-hua.fan@u-picardie.fr}

\author{J\"{O}RG SCHMELING}
\address{Lund University\\
	Centre for Mathematical Sciences\\
	Box 118, 221 00 LUND,  Sweden}
\email{joerg@maths.lth.se}

\author{Weixiao SHEN}
\address{Shanghai Center for Mathematical Sciences\\ Fudan University\\ 220 Handan Road, Shanghai 200433, China}
\email{wxshen@fudan.edu.cn}


\maketitle

\vspace{1em}

\centerline{\em Dedicated to the memory of Professor Ka-Sing Lau}

\begin{abstract} We consider the generalized Thue-Morse sequences 
	$(t_n^{(c)})_{n\ge 0}$ ($c \in [0,1)$ being a parameter) 
	 defined by 
	$
	t_n^{(c)} = e^{2\pi i c s_2(n)}
	$, where $s_2(n)$ is the sum of  digits of the binary expansion of $n$. For the polynomials 
	$\sigma_{N}^{(c)} (x) := \sum_{n=0}^{N-1} t_n^{(c)} e^{2\pi i n x}$,  we have  proved in \cite{FSS2021}  that
	the uniform norm $\|\sigma_N^{(c)}\|_\infty$  
	behaves like $N^{\gamma(c)}$ and the best exponent $\gamma(c)$ is computed. 
	In this paper, we study the pointwise behavior and give a complete multifractal analysis of the limit
	$\lim_{n\to\infty}n^{-1}\log |\sigma_{2^n}^{(c)}(x)|$.

\end{abstract}

\tableofcontents

\section {Introduction and main results}

	The first results on the multifractal formalism go already back to Besicovitch \cite{Besicovitch1934} and its systematic study has been initiated in \cite{CLP1987}, based on a physical idea introduced in \cite{HJKPS1986}. The first works are mainly concentrated on the study of local dimensions of invariant measures.
	Non-invariant measures (even without dynamics) have also been studied \cite{BMP1992,F1997}. 
	The 
	spectra for Lyapunov exponents for conformal uniformly expanding repellers have been
	covered for the first time in \cite{BPS1997} (see \cite{Pesin1997}
	for more details and references). 
	The multifractal analysis of Birkhoff averages for continuous functions were well studied for classes of nice dynamical systems including subshifts of finite type \cite{FFW2001}, 
	conformal systems \cite{FLW2002},
	systems with specifications \cite{TV2003,FLP2008}. Sub-additive potentials are also studied \cite{CFH2008}. What happens for discontinuous and unbounded functions ? This is a general problem which is worthy of study.
	The present paper contributes to this problem.  For one-dimensional 
	real and complex dynamics, multifractal analysis are performed for
	some  potential with singularities of a different kind, especially the geometric potential
	$- t \log |D f|$ \cite{MS2000,MS2003,PS2008,BT2009,PRL2011}.  

\medskip

This paper, which would have been the second part of \cite{FSS2021} and was cut off because of its length, is devoted to the multifractal analysis of the following 
ergodic averages
\begin{equation}\label{EA}
     \lim_{N\to\infty}\frac{1}{N} \sum_{n=0}^{N-1} \log |\cos \pi (2^n x +c)|   \quad (x \in [0,1))
\end{equation}
where $c\in [0,1)$ is a parameter. Let 
\begin{equation}\label{fc}
    f_c(x) := \log |\cos \pi (x +c)| = \log |\sin \pi (x-b)|  
\end{equation}
where  $ b:= \frac{1}{2}-c$.
If the function $\log |\cos \pi (x+c)|$ in (\ref{EA}) is replaced by a continuous function, the multifractal analysis of the 
ergodic average in (\ref{EA}) has been well performed (cf. \cite{FFW2001}).
But the function $f_c$ here is not continuous at $x=b$ since it has a logarithmic pole.  This singularity makes the situation more delicate.
As we shall see, the location of the singular point $b$ plays a crucial role.
 
This problem we study is motivated by arithmetic questions. 
For any integer $n\ge 0$, we denote by $s_2(n)$  the sum of  digits of the binary expansion of $n$. Fix $c\in [0,1)$, we define the generalized Thue-Morse sequence $(t_n^{(c)})_{n\ge 0}$ by
$$
        t_n^{(c)} = e^{2\pi i c s_2(n)}.
$$
The case $c=1/2$ corresponds to the classical Thue-Morse sequence
$((-1)^{s_2(n)})$.
 Then we define the {\em generalized Thue-Morse  trigonometric
polynomials} by
\begin{equation}\label{GTM}
     \sigma_{N}^{(c)} (x) := \sum_{n=0}^{N-1} t_n^{(c)} e^{2\pi i n x} \quad (N\ge 1).
\end{equation}

The Thue-Morse sequence $(-1)^{s_2(n)}$ and the  digital sum function $n \mapsto s_2(n)$, as well as the norms $\|\sigma_N^{(c)}\|_\infty$ and $\|\sigma_N^{(c)}\|_p$,  are extensively studied in harmonic analysis and number theory
after the works of Mahler \cite{M1927}  and Gelfond \cite{Gelfond1968}. Let us just cite 
\cite{AHL2017, DT2005,FK2018, FM1996a,MR2010} and the survey \cite{Queffelec2018}.
\medskip

 It is fundamental to observe 
  the following expression
 \begin{equation}
     |\sigma_{2^n}^{(c)}(x)| = 2^n \prod_{k=0}^{n-1} |\cos \pi (2^k x +c)|.
 \end{equation} 
 Thus the dynamical system  $T: \mathbb{T}\to \mathbb{T}$ defined by $Tx = 2 x \mod 1$ is naturally involved. 
 Let us define the dynamical maximum
\begin{equation}\label{Beta}
       \beta(c): = \sup_{\mu \in \mathcal{M}_T} \int_{\mathbb{T}}
         \log |\cos \pi (x+c)| d\mu(x)
\end{equation}
where $\mathcal{M}_T$ is the set of $T$-invariant Borel probability measures.  

Concerning the maximal value $\beta(c)$ and the norm $\|\sigma_N^{(c)}\|_\infty$, the following results were proved in 
\cite{FSS2021}.
\medskip

{\em 
	\indent {\rm  (I)} The supremum in (\ref{Beta}) defining $\beta(c)$
is attained by a unique measure and this measure is Sturmian. \\
\indent {\rm  (II)} Such a  Sturmian measure is  periodic
in most cases. More precisely, those parameters $c$ corresponding to non-periodic
Sturmian measures form a set of zero Hausdorff dimension. \\
\indent {\rm  (III)} There is a constant $C=C_c>0$ such that
\begin{equation}\label{eq:max}
	\forall x \in \mathbb{T},  \forall n\ge 1, \ \ \ 
	2^n\prod_{k=0}^{n-1} |\cos \pi (2^k x +c)|\le C 2^{ \gamma(c) n }
\end{equation}
where $\gamma(c): = 1 +\frac{\beta(c)}{\log 2}$ is the best possible.
}
\medskip

A dual quantity is the dynamical minimum
\begin{equation}\label{Alpha}
\alpha(c): = \inf_{\mu \in \mathcal{M}_T} \int_{\mathbb{T}}
\log |\cos \pi (x+c)| d\mu(x).
\end{equation}
Both $\alpha(c)$ and $\beta(c)$
 will play  important roles in the multifractal analysis 
of $\sigma_N^{(c)}(x)$. 
As we shall see, the minimum  $\alpha(c)$ depends strongly on how the singular 
point $b:=1/2 -c$ returns back to $b$. Let us define
$$
m_c^* := \liminf_{n\to \infty} \frac{1}{n} \sum_{k=1}^{n}\log |\sin \pi (2^k b - b)|, 
$$ 
which reflects the close returns of the orbit of $b$. The number  $m_c^* $ is closely related to $\alpha(c)$.

Concerning the minimum $\alpha(c)$, we shall prove the following theorem.
As we shall see, the set $\{c\in [0,1): \alpha(c)=-\infty\}$
is topologically large (of second Baire category), but metrically small (of zero Hausdorff dimension). 
\medskip
	
\noindent {\bf Theorem A.} {\em 
		The following assertions hold:\\
	\indent{\rm  (1)} For any $c \in [0,1)$, we have $\alpha(c) =-\infty$ iff $m_c^* =-\infty$. \\
	\indent {\rm  (2)} If $1/2 -c$ is a periodic point, then $\alpha(c)= -\infty$.\\
	\indent {\rm (3)} If $1/2 -c$ is pre-periodic but not periodic, then $\alpha(c)> -\infty$.\\
	\indent {\rm  (4)} The set of  $c \in [0,1)$ such that $
	\alpha(c) = -\infty$ is a residual set.\\
	\indent {\rm  (5)} For  Lebesgue almost all  $c \in [0,1)$, we have $ m_c^* = -\log 2$.
	\\
	\indent {\rm  (6)} $\dim_H \{c:  \alpha(c) = - \infty\} =0$. 
}
\medskip

The statements (2) and (3) follows immediately from (1). We shall prove (1) and (4) in Section 2. 
The proof of (4) is relatively easy, because the set contains all periodic points and we can 
construct a residual subset by shadowing and glueing periodic points. 
The proof of (1) uses the binding argument and the notion of first free return of   Benedicks and Carleson \cite{BC1985} in their study of the dynamics of quadratic
polynomials.  Indeed we could compare the above quantity $m_c^*$ with the Lyapunov exponent at the critical point of the quadratic polynomial. 
The proof of (5), which will be presented in Section 3,  is based on Rademacher-Menshov theorem on the almost everywhere convergence of quasi-orthogonal series (see \cite{KSZ1948}). The proof
of (6) given in Section 4 is based on the same idea used in the proof of (1) and covers used for  estimating Hausdorff dimension are carefully constructed.
\medskip

Now let us investigate the pointwise behavior  of $\sigma_{2^n}^{(c)}(x)$.
For $\alpha \in \R$ we define the level set
\begin{equation}\label{Ea}
E(\alpha):= E_c(\alpha):=\left\{ x \in \mathbb{T}: \lim_{n\to\infty} \frac{1}{n} \sum_{k=0}^{n-1}\log |\cos \pi (2^k x+c)| =\alpha \right\} . 
\end{equation}
We define the hyperbolic part of $E_c(\alpha)$ as follows
$$
 E^{\hyp}(\alpha): = E_c^{\hyp}(\alpha):=\left\{x\in E_c(\alpha): \inf_{1\le n<\infty} \dist(T^n(x), b)>0\right\}
 $$
where $d(t,s)=\|t-s\|$ with $\|x\|=\inf_{n\in \mathbb{Z}} |x-n|$  
is the usual distance  on $\mathbb{T}$.
As we shall prove, the Hausdorff dimensions of the sets $E_c(\alpha)$ are  determined by  the following pressure function
$$p(t):= p_c(t):=\sup\left\{h_\mu+ t \int f_c d\mu: \mu\in \mathcal{M}_T, \supp(\mu)\not\ni b \right\},$$
or more precisely by its Legendre-Fenchel transform
$$
     p_c^*(\alpha) := \inf_{t \in \mathbb{R}} (p_c(t) - \alpha t).
$$
Here $h_\mu$ denotes the metric entropy of $\mu$.

\medskip
Our second main theorem is as follows.


\noindent {\bf Theorem B.} {\em 
	Suppose $0<c<1$.
	The following assertions holds:\\
	 \indent{\rm  (1)}
	For $\alpha \in (\alpha(c), \beta(c))$, we have
	$$
	\dim_H E_c(\alpha) 
	= \dim_H E_c(\alpha)^{\rm hyp} =\frac{p_c^*(\alpha)}{\log 2}.
	$$
	 \indent{\rm  (2)}
	If $\alpha(c) =-\infty$, $ \dim E_c(\alpha) =1$ for all $\alpha \le -\log 2$. \\
	 \indent{\rm  (3)}
	  $\dim_H E_c(\beta(c))=0$.\\
	    \indent{\rm  (4)}
	   $ E_c(\alpha) =\emptyset$
    	if $\alpha >\beta(c)$; $\dim_H E_c(\alpha) =0$ if $\alpha < \alpha(c)$.
}
\medskip

The strategy of the proof is as follows. 
Adapting an idea of Przytycki and Rivera-Letelier \cite{PRL2007}  (in thermodynamical theory of one-dimensional maps), we decompose the set
$\{x:\omega(x)\ni b\}$ into two subsets:
 a set of $x$ with at most finitely many `good returns', which has Hausdorff dimension zero, due to the expanding property of $T$ (this set has no contribution to dimension), and  another set of $x$ with infinitely many `good returns', whose  Hausdorff dimension is bounded from above by the Legendre-Fenchel transform of the  pressure function $p(t)$ (the dimension of this set is well controlled). The complementary set of $\{x:\omega(x)\ni b\}$
 in $E_c(\alpha)$
is the hyperbolic part and can be approximated by 
irreducible and aperiodic Markov subsystems and the classic thermodynamical formalism applies to these subsets.
\medskip

Let us point out that Theorem B fails for the case $c=0$ (i.e. $b = 1/2$). This particular case will be discussed in Section 5. When $c=0$, we have 
$$
\forall t \in \mathbb{R}, \ \ \ p_0(t) = \max \{(1- t)\log 2,0\};
$$ 
$$
p_0^*(\alpha) = |\alpha|
\ \  {\rm for } \   \alpha \in [-\log 2, 0];
\quad p_0^*(\alpha) = -\infty\  {\rm for \ other } \ \  \alpha.
$$
  The Birkhoff limits avoid the interval $
(\alpha(0), \beta(0))$, 
i.e. $E_0(\alpha)=\emptyset$
for $\alpha \in (-\log 2, 0)$. Actually as we shall see,   there are only three possible limit values 
$-\infty, -\log 2$ and $0$.
The proof of Theorem B will be given in Section 6. 
\medskip

Finally we point out that our works can easily be adopted  to generalized Thue-Morse trigonometric polynomials defined by the
digital sum function $s_q(n)$ in base $q\ge 2$.

\medskip

{\bf Acknowledgements.} 
The first author is supported by NSFC grant no. 11971192 and the
third author is supported by NSFC grant no. 11731003. The first and second authors would like to thank  
Knuth and Alice Wallenberg Foundation (2014-2016) and  Institut Mittag-Leffler (Oct. 2017)
 for their supports.

\section{Minimization for $f_c$ and exponent $\chi_-(f_c, T(b(c)))$}
Recall that $f_c(x) = \log |\cos \pi (x+c)| = \log |\sin \pi (x-b)|$  
with $b:= b(c):=\frac{1}{2}-c$. The $1$-periodic function $f_c$ is real analytic except at the 
singular point $b(c)$ where we have $f_c(b(c))= - \infty$.
Recall that  $$\alpha (c)=\inf\left\{\int f_c d\mu: \mu\in \mathcal{M}_T\right\}.$$ It is not always possible to show that a minimizing measure exists.
The maximum $\beta(c)$ was well studied din \cite{FSS2021}. In this section we would like to determine those parameters $c$ such that 
$
\alpha(c) = -\infty.
$
The discussion is very similar to the study of Collet-Eckmann maps in interval dynamics and the singular point $b(c)$ plays the role of the critical point (see \cite{CE1980,BC1985}). 

For each $x\in \T$, define the exponent
$$\chi_-(f_c,x)=\liminf_{n\to\infty} \frac{1}{n}\sum_{i=0}^{n-1} f_c(T^i(x)),$$
and $\chi_+(f_c,x)$ similarly. If $\chi_-(f_c,x) =\chi_+(f_c,x)$,
the common value is denoted by $\chi(f_c,x)$.
These quantities have obvious similarity with Liapunov exponents.
As we shall show, $\alpha(c) =-\infty$ iff $\chi_{-}(f_c, T(b(c)))=-\infty$ and that
the set of $c$'s such that $\alpha(c) = -\infty$ is uncountable but has zero Hausdorff dimension. 

Remark that instead of $f_c$, we can work with $\log_2 d(x, b(c))$ if we are only interested in the finiteness of $\chi_-(f_c, T(b(c)))$, where 
$d(x, y)$ is the distance on $\mathbb{T}$ define by $d(x, y) = \|x-y\|$.
This is 
because 
$$
\forall x\in [-1/2, 1/2], \ \ \ \frac{2}{\pi}\|x\|\le |\sin \pi x|\le \pi \|x\|.
$$

In Section \ref{sect:Existence}, we shall also prove that for Lebesgue almost all $c$ we have
$$
    \chi_-(f_c, T(b(c))) = \chi_+(f_c, T(b(c))) = \int f_c (x) dx =- \log 2.
$$

We first prove that both $\alpha(c)$ and $\beta(c)$ are determined by periodic invariant measures.
  
\subsection{ $\alpha(c)$ and $\beta(c)$}

\begin{prop} \label{prop:ab}
Let $\mathcal{P}_T$ be the set of all $T$-invariant probability measures supported on periodic orbits. We have
$$
    \alpha(c) = \inf_{\mu \in \mathcal{P}_T} \int f_cd\mu, \qquad \beta(c) = \sup_{\mu \in \mathcal{P}_T} \int f_c d\mu.
$$
\end{prop}
\begin{proof} The result on $\beta(c)$ was already proved in \cite{FSS2021} (cf. Theorem 2.1 there).
Let us prove the result on $\alpha(c)$. For simplicity, we write $f$ for $f_c$.

	For each positive integer $N$, consider the continuous function $f_N(x)=\max\{f(x), - N\}$, which decreases to $f$ when $N \uparrow \infty$. 
 	Then consider the closed set $G_N =\{x: f(x) \ge -N\}$ and let  $$\mathcal{F}_N = \bigcap_{n=0}^\infty T^{-n} G_N.$$ 
 	Clearly $\mathcal{F}_N$ is non-empty and closed,  and $T(\mathcal{F}_N)\subset \mathcal{F}_N$.
 	Let 
 	$\mathcal{P}(\mathcal{F}_N)$ be the set of invariant periodic probability measures of the subsystem $T: \mathcal{F}_N \to \mathcal{F}_N$. 
 	Let 
 	$$
 	\alpha_{f, N}(f) := \inf_{\mu \in \mathcal{P}(\mathcal{F}_N)} \int f_N d\mu.
 	$$
 	Since $ \mathcal{P}(\mathcal{F}_N)\subset  \mathcal{P}_T$
 	and  $f_N \ge f$, we get immediately 
 	\begin{equation}\label{MM1}
 	\alpha_f \le \inf_{\mu \in \mathcal{P}_T} \int f d\mu \le \inf_N \alpha_{f, N}.
 	\end{equation}  
 	We distinguish two cases for proving the reverse inequality of (\ref{MM1}). 
 	
 	{\em Case 1:  The singular point $b$ is a periodic point}. Notice that $b$ is not a periodic point of any subsystem $\mathcal{F}_N$, because $b \not \in \mathcal{F}_N$. But we can construct 
 	a  periodic point for the subsystem, which returns often close to $b$. Let $b = \overline{b_1b_2\cdots b_\tau}$ where $\tau$ is the  period of $b$. Let 
 	$$
 	b^* = \overline{(b_1\cdots b_\tau)^{N-1} b_1\cdots b_{\tau-1}\check{b}_{\tau}},
 	$$
	where $\check{b}_{\tau}=0$ or $1$ according to $b_\tau=1$ or $0$.
  Then $b^*$  is $\tau N$-periodic and $b^* \in \mathcal{F}_{\tau N}$ because
 	$$
 	     f(b^*) \ge c \log |b^* - b| \ge - (\tau N+1) c\log 2.  
 	$$
 	Let $\mu$ be the periodic measure supported by the orbit of  $b^*$. Then we have
 	$$
 	\int f_{\tau N} d\mu \le \frac{1}{\tau N}\sum_{j=0}^{N-1} f(T^{j \tau} b^*).
 	$$
 	Notice that $T^{j\tau}b^*$ is close to $b$ at least for $j$ not too large.  Then the right hand side goes to $-\infty$ and we get 
 	$$
 	 \inf_N \alpha_{f,N} = -\infty =\alpha_f.
 	$$

 	{\em Case 2: The singular point $b$ is not periodic.}  
 	Assume $\alpha_f >-\infty$ (the case $\alpha_f =-\infty$ can be   similarly proved). For any $\epsilon>0$, there exists $\mu_0 \in \mathcal{M}_T$ such that $\alpha_f> \int f d\mu_0 - \epsilon$. 
 	Take a large integer $N_0\ge 1$ such that  $\int f d\mu_0\ge \int f_{N_0} d\mu_0 - \epsilon$. Since $f_{N_0}$ is continuous,  we can find a periodic measure $\nu_0 \in \mathcal{P}_T$ such that $\int f_{N_0} d\mu_0 >\int f_{N_0} d\nu_0 -\epsilon$. Hence we have  
 	$$
 	\alpha_f > \int f_{N_0} d\nu_0 -3 \epsilon 
 	\ge \int f_{N} d\nu_0 -3 \epsilon \quad (\forall N \ge N_0). 
 	$$
 	The support of $\nu_0$ must have a positive distance from the singular point $b$, otherwise  $\int f_{N} d\nu_0$ tends to $-\infty$. So,  $\nu_0 \in  \mathcal{M}(\mathcal{F}_N)$ for $N$ sufficiently large. Hence, for large $N$  we have
 	$$
 	\alpha_{f,N} \le \int f_N d\nu_0 <\alpha_f +3\epsilon.
 	$$
 	Thus $\inf_N \alpha_{f, N} \le \alpha_f$. 
 \end{proof}

\subsection{ $\{c: \alpha(c) =-\infty\}$ is a residual set}
Assuming that Theorem A (1) is proved, 
 we give here a proof of Theorem A (4), in the symbolic setting. The set in question contains all periodic points. By gluing periodic points, we construct a subset which is residual. Let $c=(c_k)_{k\ge 1}$ be a $n$-periodic point.
 Define  the open set
 $$
 A_n(c)=\{x \in \{0,1\}^\infty : x_1x_2\cdots x_{n^2} = (c_1...c_n)^n\}.
 $$
 The following set is a residual set
 $$
 A_\infty:=\bigcap_{N\ge 1} \ \ \ \bigcup_{c \ {\rm is} \ n\!-\!{\rm  periodic}, n\ge N} A_n(c).
 $$ 
 We claim that $\alpha(x) =-\infty$ for all $x\in A_\infty$.
  
 Notice that if $d(x,y)\le 2^{-k}$, 
  we have 
     $
           \log |\sin \pi (y-x)|\le - C k
     $
 where $C>0$ is an absolute constant. 
It follows that for $x \in A_n(c)$, we have 
\begin{equation}\label{Residue}
     S_{n^2}f_x(Tx)\le - \frac{C}{2} n^3 
\end{equation}
because $d(2^{nj} x, x)\le 2^{-{n^2-jn}}$ and 
$$
S_{n^2}f_x(Tx) \le \sum_{j=1}^{n-1} \log |\sin \pi (2^{jn}x -x)| 
\le - C \sum_{j=1}^{n-1} (n^2 -jn). 
$$ 
Here we have just kept $n-1$ main terms in the sum defining $S_{n^2}f_x(Tx)$.
If $ x\in A_\infty$, (\ref{Residue}) holds for infinitely many 
$n$'s. So $m^*_x = -\infty$.

\subsection{$\alpha(c) =-\infty$ is equivalent to $\chi_-(f_c,T(b(c)))=-\infty$}\
In other words, we prove here the statement (1) of Theorem A. 
It is clear that if $b(c)$ is a periodic point of $T$, then $\alpha (c)=-\infty$. More generally, we have
\begin{prop}\label{chi_infty}
If $\chi_-(f_c, T(b(c)))=-\infty$, 
then $\alpha(c)=-\infty.$
\end{prop}
\begin{proof} 
Let $q_n:=-\log_2 \dist (T^n(b(c)), b(c))$. 
Since $$
\log |\sin \pi x| = \log \|x\| + O(1)\quad {\rm  when }\ \ \ |x|\le \frac{1}{2},
$$ we have $f_c(T^{n}(b(c)))=-Aq_n +O(1)$ with $A=\log 2$, so that
  \begin{equation}\label{-infty}
  \chi_-(f_c, T(b(c)))=-\infty  \Leftrightarrow \limsup_{n\to\infty} \frac{1}{n}\sum_{i=1}^n q_i =\infty.
\end{equation}

We distinguish two cases.

{\bf Case 1.} $\sup \frac{q_n}{n}=\infty$. This is the case when $b$ can return very closely to $b$. 
Let $n_1<n_2<\cdots$ be integers such that $\frac{q_{n_k}}{n_k}\to\infty$. For each $k$, let $j_k\in \Z$ be such that
$$|(2^{n_k}-1)b -j_k| =|2^{n_k}b-b-j_k|\asymp  2^{-q_{n_k}}$$
and let $p_k=j_k/(2^{n_k}-1)\mod 1$, which is a $n_k$-periodic point of $T$,  close to $b$ in the sense that $\dist(b, p_k)\le 2^{-q_{n_k}}$. It follows that 
$$\chi_{-}(f_c, p_k)
\le \frac{ f_c(p_k)}{n_k} \le -A \ \frac{q_{n_k}}{n_k} +O(1)
\to -\infty.$$  That is to say  $\int f_c d\mu_k\to-\infty$, where  $\mu_k$ is the invariant measure supported by the orbit of $p_k$.

{\bf Case 2.} $Q:=\sup \frac{q_n}{n} <\infty$. Let $m>Q+1$ be a large integer. For any constant $C>0$, we shall find a $T$-invariant probability measure $\mu$ such that $\int \log_2 \dist (x, b(c)) d\mu<-C$, which will complete the proof. 
To this end, let $N$ be a large enough integer so that
$$\frac{1}{N}\sum_{i=1}^N q_i> Cm+1.$$
Such an $N$ exists according to (\ref{-infty}).  
 Let $I=[j/2^{mN}, (j+1)/2^{mN})$ be the dyadic interval containing $b(c)$, where $j\in \{0,1,\cdots, 2^{mN}-1\}$ and let $p$ be the $mN$-periodic  point of $T$ contained in the closure of $I$. 
 Then $|p-b(c)|\le 2^{-mN}$ and even
  $T^k(p)$ follows closely $T^k(b)$ when $1\le k\le N$: indeed, from the fact
 $$|T^k(p)-T^k(b(c))|\le 2^{k} |p-b(c)|\le 2^{-(m-1)N},$$
 and the fact $m-1 >Q$, we have
 $$ |T^k(p)-T^k(b(c))| \le 2^{-Q N}\le 2^{-Q k}\le |T^k(b(c))-b(c)|.$$
 Then, using the triangle inequality, we get 
 $$|T^k(p)-b(c)|\le 2 |T^k(b(c))-b(c)|.$$
Thus 
$$\sum_{k=1}^N \log_2 \dist (T^k(p), b(c))\le N -\sum_{k=1}^N q_k< -CmN.$$
Hence 
$$\sum_{k=1}^{mN} \log_2 \dist (T^k(p), b(c))\le \sum_{k=1}^N \log_2 \dist (T^k(p), b(c))\le -CmN.$$
Let $\mu$ be the $T$-invariant probability measure supported by the
$mN$-periodic orbit of $p$. What we have proved means 
$$\int \log_2 d(x, b(c)) d\mu =\frac{1}{mN} \sum_{k=1}^{mN} \log_2 \dist (T^k(p), b(c))\le -C.$$  
\end{proof}

The converse of Proposition \ref{chi_infty} is also true. 
\begin{prop}\label{prop:lbarb} Suppose $\chi_-(f_c,T(b(c))>-\infty$. Then
	there is a constant $K>0$ such that for every $x\in \T$
	either  $T^n(x)=b(c)$ for some $n\ge 0$, or $\chi_+(f_c, x) >-K$.
Consequently, $\alpha(c) \ge -K$. 
\end{prop}
The following argument in our proof of Proposition \ref{prop:lbarb} is inspired by the `binding argument', which is   developed by Benedicks-Carleson \cite{BC1985} to study iterations of Collet-Eckmann\cite{CE1980} quadratic polynomials and which has been extremely influential in the study of iterations of non-uniformly hyperbolic maps.
\medskip

The assumption implies that $b=b(c)$ is not periodic under $T$ and that there exists a constant $K_0>0$  such that
\begin{equation}\label{bind0}
\forall n \ge 1, \quad \frac{1}{n}\sum_{i=1}^{n} \log_2 \dist (T^n(b), b)>-K_0.
\end{equation}
For each $\rho>0$, let $p:=p(\rho)$ be the maximal non-negative integer such that
\begin{equation}\label{bind1}
\forall\ 1\le k\le p, \quad  2^{k+1}\rho < \dist (T^k(b), b).
\end{equation}
The maximality implies
\begin{equation}\label{bind2}
 2^{p+2} \rho \ge d(T^{p+1}b, b). 
 \end{equation}
 
 If $\rho \ge \frac{1}{4}d(Tb, b)$, we have $p(\rho)=0$. But $p(\rho)\ge 1$
 if $\rho < \frac{1}{4} d(Tb, b)$. 
It is clear that
there exist $\rho_0: =\frac{1}{8}d(Tb, b)>0$ and $C:= \frac{1}{2(K_0+2))}>0$ such that for each $\rho\in (0, \rho_0]$,
\begin{equation}\label{eqn:prho}
p(\rho)\ge C \log\frac{1}{\rho}.
\end{equation}
This follows from (\ref{bind0}) applied to $n=p+1$ and (\ref{bind2}):
$$
  -K_0 < \frac{\log_2 d(T^{p+1} b, b)}{p+1} \le \frac{p+2 + \log_2 \rho}{p+1}\le 2+ \frac{\log_2 \rho}{p+1}.
$$

A  reversed inequality of (\ref{eqn:prho}) is also true, because of (\ref{bind1}) from which we get
$$\sum_1^p \log_2 (2^{k+1}\rho) < \sum_1^p \log_2 d(T^k b, b)\le 0.$$
So    $\sum_1^p (k+1) + p \log_2 \rho \le 0$,
i.e.    $  (p+2)(p-1)/2 + p \log_2 \rho \le 0$.

\begin{definition} For each $x\not=b$, we say that $p(\dist(x,b))$ is {\em the binding period of $x$}. The minimal positive integer $s>p(\dist(x,b))$ with $\dist(T^sx, b)<\rho_0$ is called the first {\em free return of $x$} (if it exists).
\end{definition}

\begin{lemma} \label{lem-bind}
	For each $x\not=b$, if $s$ is the first free return time of $x$, then
$$\frac{1}{s}\sum_{i=0}^{s-1} \log_2 \dist (T^i(x), b)\ge -K_1,$$
where $K_1$ is a constant independent of $x$.
\end{lemma}
\begin{proof} Let $\rho=\dist(x,b)\not=0$ and $p=p(\rho)$. By the definition of $s$, we have $\log \dist(T^i(x), b)\ge \log \rho_0$ for each $p(\rho)<i<s$. So it suffices to show that there is a constant $K_1'>0$ such that
\begin{equation}\label{eqn:bindingest}
\sum_{i=0}^{p} \log_2 \dist (T^i(x), b)\ge -K_1' (p+1).
\end{equation}
We may assume $\rho<\rho_0$. Otherwise, let $q$ be the smallest integer such that
$d(T^q x, b) <\rho_0$. Then
$$\sum_0^{q-1} \log_2 d (T^i x, b) \ge q \log_2 \rho_0$$
and we only need to discuss $T^q x$ instead of $x$.

Then for each $1\le k\le p$,
\begin{align*}
\dist (T^k(x), b)& \ge \dist (T^k (b), b)-\dist(T^k(x), T^k(b))\\
& \ge\dist(T^k(b), b)-2^k \dist(x,b)\ge \dist(T^k(b), b)/2.
\end{align*}
Here for the last inequality, we have used (\ref{bind1}). Hence

$$\log_2 \dist (T^k(x), b)\ge \log_2 \dist(T^k(b), b)-1.$$ Thus
\begin{align*}
\sum_{i=0}^{p} \log_2 \dist (T^i(x), b)& =\log_2 \rho +\sum_{i=1}^{p}\log_2\dist(T^i(x), b)\\
& \ge \log_2\rho +\sum_{i=1}^{p} \left(\log_2 \dist(T^i b, b)-1\right)\\
& \ge  \log_2\rho -K_0p - p.
\end{align*}
Since $p\le C\log_2\rho^{-1}$, the estimate (\ref{eqn:bindingest}) follows.
\end{proof}
\begin{proof}[Proof of Proposition~\ref{prop:lbarb}]
Let $x\in \T$ be such that $T^n x\not=b$ for all $n\ge 0$. We may assume that there are infinitely many $n$ such that $\dist(T^n(x), b)<\rho_0$, for otherwise, 
$\dist(T^n(x), b)\ge \rho_0$ for all $n\ge n_0$ (some $n_0$) and 
$$\chi_+(f_c, x)\ge \min\{f_c(T^n(x)): n\ge n_0\}=:-K_2>-\infty.$$
We define inductively a sequence of integers $0\le s_0<s_1<\cdots$ as follows:
\begin{itemize}
\item $s_0$ is the first free return time of $x$;
\item for each $k\ge 0$, $s_{k+1}$ is the first free return time of $T^{s_k}(x)$.
\end{itemize}
Then by Lemma \ref{lem-bind},
$$\sum_{i=s_k}^{s_{k+1}-1} \log_2\dist (T^i(x), b) \ge -K_1 (s_{k+1}-s_k).$$
It follows that $\chi_+(f_c,x)\ge -K_1$. Take $K=\max(K_1,K_2)$. Then 
$\chi_+(f_c,x)\ge -K$ holds for all $x$ such that $T^n x \not= b$
for all $n\ge 0$. For any (ergodic) periodic measure $\mu$, by Birkhoff
ergodic theorem, we have
$$
    \int f_c(x) d\mu = \chi_+(f_c, x) \ \ \  \mu\!-\!a.e.
$$
As $\chi_+(f_c, x) \ge -K$ $\mu$-a.e., we get 
 $$
 \alpha(c) = \inf_{\mu \ {\rm periodic}}   \int f_c(x) d\mu \ge -K,$$
 where the equality is ensured by Proposition \ref{prop:ab}.
\end{proof}

\section{Existence of the exponent $\chi (f_c, T(b(c)))$}\label{sect:Existence}
We are going to prove the statement (5) of Theorem A: relative to the Lebesgue measure 
we have
\begin{equation}\label{eq:ExpChi}
   a.e. \ \ \  \chi(f_c, T(b(c))) = - \log 2. 
\end{equation}
Actually we shall prove a more general result. 
For $f \in L^1(\mathbb{T})$, recall that the $L^1$-modulus of continuity of $f$ is defined by 
$$
\Omega_1(f; \delta) := \int_{\mathbb{T}}|f(x+\delta) - f(x)| d x.
$$

\begin{theorem} 
	For any $f \in L^2(\mathbb{T})$ such that $\Omega_1(f) = O(\delta^{1/2 +\eta})$ for some $\eta>0$ we have
	$$
	a.e. \quad   
	\lim_{N\to \infty} \frac{1}{N} \sum_{n=1}^N f(2^n \xi - \xi) = \int f(x) d x. 
	$$
\end{theorem}

\begin{proof} We can not apply the Birkhoff ergodic theorem. But we can apply the Rademacher-Menshov theorem on quasi-orthogonal systems \cite{KSZ1948}.  
Since the estimate $|\widehat{f}(n)|\le \Omega_1(f, \pi/|n|)$ holds for every integrable function $f$,  for our function $f$ we have
	\begin{equation}\label{OmEstimate}
	\widehat{f}(n) = O(n^{-1/2-\eta}).
	\end{equation}
	for any integer $n \not=0$.
We can assume 
$\int f(x) dx =0$ without loss of generality
	and 
	consider the covariance function of $f$ defined by 
	$$
	\mbox{\rm Cov}(f; p, q)  
	:=   \int_0^1  f(px) \overline{f(qx)} dx
	$$
	for integers $p$ and $q$. We have only to show that  $$
	\mbox{\rm Cov}(f; 2^{j}-1, 2^{k}-1) = O(|k-j|^{-2}),$$
	because this implies that $f((2^n-1)\xi)$ ($n =1, 2, \cdots$) is a quasi-orthogonal system (cf. \cite{KSZ1948}) to ensure 
	the almost everywhere convergence of the series $\sum n^{-1} f(2^n \xi - \xi)$. This implies the desired result via the Kronecker lemma.  
	
	We shall prove that $\mbox{\rm Cov}(f; 2^{j}-1, 2^{k}-1)$
	has an exponential decay as function of $|k-j|$.   
	Denote $f_p(x) = f(px)$. By Plancherel's formula we have
	\begin{equation*}\label{OmEstimate0}
	\mbox{\rm Cov}(f; p, q)  = \sum_{n\not =0} \widehat{f}_p(n) \overline{\widehat{f}_q(n)}.
	\end{equation*}
	Notice that 
	$\widehat{f}_p(n) = \widehat{f}(n/p)$ if $p | n$ and 
	$\widehat{f}_p(n) =0 $ if $p \not| n$. 
	So, if  $[p,q]$ denotes the common least multiple of $p$ and $q$, we have
	\begin{equation}\label{OmEstimate2}
	\mbox{\rm Cov}(f; p, q)  =  \sum_{m\not =0} \widehat{f}([p,q]m/p) \overline{\widehat{f}([p,q]m/q)}.
	\end{equation}
	Thus, from  (\ref{OmEstimate2}) and  (\ref{OmEstimate}),
	we have
	\begin{eqnarray*}
		| \mbox{\rm Cov}(f; p, q)| & \le &    \sum_{m\not =0} \Omega_1\left(f; \frac{\pi p}{[p,q]|m|} \right) \Omega_1\left(f; \frac{\pi q}{[p,q]|m|} \right)\\
		& \le &   O\left( \left(\frac{pq}{[p,q]^2}\right)^{1/2+\eta}  \sum_{m=1}^\infty \frac{1}{m^{1+ 2\eta}}\right) =  O\left( \left(\frac{pq}{[p,q]^2} \right)^{\frac{1}{2}+\eta}\right).
	\end{eqnarray*}
Let $p=2^j -1$ and $q=2^{k} -1$ with $1\le j<k$. Then 
$$\frac{pq}{[p,q]^2}\le \frac{p}{q}\le 2/2^{k-j}.$$ So  we get the exponential decay 
	$$
	\mbox{\rm Cov}(f; 2^j -1, 2^k -1)  =O\left(  2^{-(1/2+\eta) |k-j|}\right).
	$$

	

\end{proof}

We can apply the above theorem to $f(x) = \log |\sin \pi x|$
to get (\ref{eq:ExpChi}) by using the following result.

\begin{theorem} \label{thm:finiteA} 
	There exist constants $0<C_1<C_2<\infty$ such that
	\begin{equation}\label{Omega}
	C_1 \delta |\log \delta| \le \Omega_1(\log |\sin \pi x|; \delta) \le C_2 \delta |\log \delta|.
	\end{equation}
	Consequently 
	for almost all $\xi$ with respect to Lebesgue measure, we have
	\begin{equation}\label{Omega2}
	\lim_{N\to \infty} \frac{1}{N} \sum_{n=1}^N \log |\sin \pi (2^n \xi-\xi)|= - \log 2. 
	\end{equation}
\end{theorem}
\begin{proof} 
	Let us integrate on the interval $[-\delta, 1-\delta]$ and  consider the decomposition
	$$\Omega_1(f;\delta)=\left(\int_{-\delta}^0 +\int_0^{1/2-\delta}+\int_{1/2-\delta}^{1/2} +\int_{1/2}^{1-\delta}\right)|f(x+\delta)-f(x)| dx
	$$
	where $f(x) = \log |\sin \pi x|$, and then estimate each of these four integrals.
	Let  $F(x)=\int_0^x f(t) dt$. We make the following useful remark that 
	$$F(\delta)=F(1)-F(1-\delta) \asymp  \delta |\log \delta|, 
	\quad F(1/2)-F(1/2-\delta) = O(\delta)$$
	where  the notation $F(\delta) \asymp G(\delta)$ means that there exist two constants 
	$0<A <B<\infty$ such that $A |G(\delta)| \le |F(\delta)|\le B |G(\delta)|$ for small $\delta >0$.

	Let us first prove  the second estimate in (\ref{Omega}), which allows us to conclude for (\ref{Omega2}). 
	Since $f(x)$ is monotone increasing in $[0,1/2]$, we have 
	\begin{align*}
	\int_0^{1/2-\delta} |f(x+\delta)-f(x)| dx & =\int_0^{1/2-\delta}(f(x+\delta)-f(x))dx\\
	& =
	F\left(\frac{1}{2}\right)-F\left(\frac{1}{2}-\delta\right)-F(\delta)\\
	& \asymp O(\delta)+ \delta |\log\delta| \asymp \delta|\log\delta|).
	\end{align*}
	Similarly, 
	$$\int_{1/2}^{1-\delta} |f(x+\delta)-f(x)|dx = O( \delta |\log \delta|).$$
	Moreover,
		\begin{align*}
	\int_{-\delta}^0 |f(x+\delta)-f(x)|dx &\le  \int_{-\delta}^0 |f(x+\delta)| dx+\int_{-\delta}^0 |f(x)| dx\\
	&= -2 \int_0^\delta f(x) dx = O( \delta |\log \delta|),
	\end{align*}
	and since $f$ is smooth in $(0,1)$, 
	$$
	\int_{1/2-\delta}^{1/2} |f(x+\delta)-f(x)| dx=O(\delta^2).
	$$ 
Therefore,
$\Omega_1(f;\delta)=O(\delta|\log\delta|)$.

To get the first estimate in (\ref{Omega}), it suffices to use 
$$
   \Omega_1(f;\delta) \ge \int_0^{1/2-\delta} |f(x+\delta)-f(x)| dx
$$
and the above estimation for the last integral.
\end{proof}

\section{Hausdorff dimension of  $\{c: \chi_-(f_c, T(b(c)))=-\infty\}$}

Let 
$$X_\infty=\left\{x\in \T: \liminf_{n\to\infty} \frac{1}{n} \sum_{i=1}^n \log_2 \dist(T^ix, x)=-\infty\right\}.$$

\begin{theorem}\label{thm:zerodim} The set $X_\infty$ has Hausdorff dimension zero.
\end{theorem}

 In other words, the set $\{c: \chi_-(f_c, T(b(c))=-\infty\}$ or equivalently the set $\{c: m_c^*=-\infty\}$ has Hausdorff dimension zero
 (cf. Theorem A (1)). 

We shall construct suitable covers of $X_\infty$ by dyadic intervals.
Let $\mathcal{D}_k$ be the collection of $k$-level dyadic intervals of the form $[m/2^k, (m+1)/2^k)$.
For $i\ge 1$ and $j\ge 1$, let $$Q_{i}^j=\{x\in \T: \dist(T^i x,x)\le 2^{-j}\}$$ and let $\widehat{Q}_i^j$ denote the union of all intervals in  $\mathcal{D}_{i+j}$ that intersect $Q_i^j$.

\begin{lemma}\label{lem:qij} 
There exists a positive constant $M:=15$ such that for any integers $i, j\ge 1$ and for any $J\in\mathcal{D}_i$, $J\cap \widehat{Q}_i^j$ is a union of at most $M$ intervals in $\mathcal{D}_{i+j}$. Consequently, $\widehat{Q}_i^j$ consists of at at most $M2^i$ intervals in $\mathcal{D}_{i+j}$ and $|Q_i^j|\le |\widehat{Q}_i^j|\le M 2^{-j}$.
\end{lemma}
\begin{proof} Let $J=[m/2^i, (m+1)/2^i)\in \mathcal{D}_i$. 
	Assume $x\in J \cap Q_i^j$. First $x\in J$ implies $0\le 2^i x - m <1$, i.e.
	$T^i x = 2^i x -m$, hence $2^i x -m -x \in (-1, 1)$. Next  $d(T^ix, x)<2^{-j}$ further implies that one of the following three cases holds
$$ 
|2^i x-x-m-\eta|<2^{-j}, \quad (\eta=0, -1, \ {\rm or}\  1).
$$ 
Thus $J\cap Q_i^j$ is contained in 3 intervals $U_k$ of length $\frac{2}{2^j(2^i-1)}$, $k=1,2,3.$ Since $|U_k|< 4\cdot 2^{-(i+j)}$, $U_k$ can intersect at most $5$ intervals in $\mathcal{D}_{i+j}$.
Consequently, the number of intervals in $\mathcal{D}_{i+j}$ which intersect  $J\cap \widehat{Q}_i^j$ is bounded from above by $M:=15$.
\end{proof}

Let us first compute the Hausdorff dimension of a subset of $X_\infty$,
which consists of points returning very closely to themselves and which is defined by 
$$Y_{\infty}=\{x\in \T: \liminf \frac{1}{n}\log_2\dist (T^n x, x)=-\infty\}.$$
To this end, we observe that $Y_\infty =\bigcap_{K=1}^\infty Y_K$ where
$$Y_K=\{x\in \T: \liminf \frac{1}{n}\log_2\dist (T^n x, x)< -K\}.$$

\begin{prop}\label{lem:dimY} We have $\dim_H Y_K\le \frac{1}{K+1}$. Hence $ \dim_H Y_\infty =0$.
\end{prop}
\begin{proof} For each $K>0$, we have
$$Y_K \subset \bigcap_{N=1}^\infty \bigcup_{n=N}^\infty Q_n^{Kn}.$$
For each $n\ge 1$, $Q_n^{Kn}$ is contained in at most $M2^n$ intervals in $\mathcal{D}_{(K+1)n}$, by Lemma \ref{lem:qij}. Thus
for $\alpha>1/(K+1)$, the fact
$$
   \sum_{n=N}^\infty M 2^n \cdot 2^{-(K+1)n \alpha} \to 0 \ \ \ {\rm as} \ \  N \to \infty
$$
implies that 
the $\alpha$-dimensional Hausdorff measure of $Y_K$ is zero. The lemma follows.
\end{proof}

For $1\le i_1<i_2<\cdots <i_k$ and $j_1, j_2,\cdots, j_k\in \Z^+$, let
$$
Q_{i_1i_2\cdots i_k}^{j_1j_2\cdots j_k}:=\bigcap_{\nu=1}^k Q_{i_k}^{j_k}=\{x\in \T: \dist (T^{i_\nu} x, x)\le 2^{-j_\nu}, \nu=1,2,\ldots, k\},
$$
and let $$\widehat{Q}_{i_1i_2\cdots i_k}^{j_1j_2\cdots j_k}:=\bigcap_{\nu=1}^k\widehat{Q}_{i_k}^{j_k}.$$

\begin{lemma}\label{lem:covering} If $i_{l+1}-i_l\ge j_l$ holds for each $1\le l<k$, then
$\widehat{Q}_{i_1i_2\cdots i_k}^{j_1j_2\cdots j_k}$ consists of at most  $M^k 2^{i_k-j_1-j_2-\cdots-j_{k-1}}$ intervals in $\mathcal{D}_{i_k+j_k}$ ($j_0=0$).
\end{lemma}
\begin{proof} We prove this lemma by induction on $k$. The case $k=1$ is ensured by  Lemma~\ref{lem:qij}. For the induction step, let $k\ge 2$ and suppose that $\widehat{Q}_{i_1i_2\cdots i_{k-1}}^{j_1j_2\cdots j_{k-1}}$ consists of at most $M^{k-1} 2^{i_{k-1}-j_1-\cdots-j_{k-2}}$ intervals in $\mathcal{D}_{i_{k-1}+j_{k-1}}$. Each of these intervals, which is of length 
	$2^{-(i_{k-1}+j_{k-1})}$,  is the union of $2^{i_k-(i_{k-1}+j_{k-1})}$ intervals in $\mathcal{D}_{i_k}$. By Lemma~\ref{lem:qij}, each interval in $\mathcal{D}_{i_k}$ intersects $\widehat{Q}_{i_k}^{j_k}$ in at most $M$ intervals in $\mathcal{D}_{i_k+j_k}$. It follows that $$\widehat{Q}_{i_1i_2\cdots i_k}^{j_1j_2\cdots j_k}=\widehat{Q}_{i_1i_2\cdots i_{k-1}}^{j_1j_2\cdots j_{k-1}}\cap \widehat{Q}_{i_k}^{j_k}$$ consists of at most 
	$$
	M^{k-1} 2^{i_{k-1}-j_1-\cdots-j_{k-2}} \cdot M2^{i_k-(i_{k-1}+j_{k-1})}=M^k 2^{i_k-j_1-j_2-\cdots -j_{k-1}}
	$$
	 intervals in $\mathcal{D}_{i_k+j_k}$.
\end{proof}

For $\eps>0$, let $\mathcal{F}_k(\eps)$ denote the collection of sequences of positive integers $(i_1i_2\cdots i_k; j_1j_2\cdots j_k)$ for which the following hold:
\begin{enumerate}
\item[(i)] for each $1\le l<k$, $i_{l+1}-i_l \ge j_l$;
\item[(ii)] $k\le \eps i_k$.
\item[(iii)] $i_k-\sum_{1\le l<k} j_l\le \eps i_k$.
\end{enumerate}
We  call $j_l$'s spacings, according to (i). (ii) means that there are not too many spacings, but big spacings by (iii).
  
Let $\mathcal{F}_k^{n,m}(\eps)$ denote those elements of $\mathcal{F}_k(\eps)$ with $i_k=n$ and $j_k=m$. Let $\mathcal{F}^{n,m}(\eps):=\bigcup_{k=1}^n \mathcal{F}_k^{n,m}(\eps)$.
We shall need the following combinatorial lemma for counting $\mathcal{F}^{n,m}(\eps)$.

\begin{lemma}\label{lem:counting} There exist $C>0$ and $\eps'>0$ such that for each $n,m\ge 1$,
$$\# \mathcal{F}^{n,m}(\eps)\le C 2^{\eps' n}.$$
Moreover, $\eps'\to 0$ as $\eps\to 0$.
\end{lemma}
\begin{proof} If $\mathcal{F}_k^{n,m}(\eps)\not=\emptyset$, then $k\le [\eps n]$, according to (ii). For each fixed $k$, the choices of $(i_1i_2\cdots i_k)$ are not more than
$\binom{n-1}{k-1}$, $i_k = n$ being fixed. 
For $(i_1i_2\cdots i_k, j_1j_2\cdots j_k)\in \mathcal{F}_k^{n,m}(\eps)$,  let 
$$s(j_1,j_2,\cdots, j_k):=\sum_{l=1}^{k-1} (i_{l+1}-i_l-j_l).$$
By (iii), we have
$$s(j_1,j_2,\cdots, j_k) =i_k-i_1-\sum_{1\le l<k} j_l \le i_k-\sum_{1\le i<k} j_l \le \eps n.$$
For fixed $1\le s\le \epsilon n$ and fixed $(i_1,j_2,\cdots, i_k)$, the choices of possible $(j_1,j_2,\cdots, j_k)$ such that $s(j_1,j_2,\cdots, j_k)=s$ are not  more than non-negative integer solutions of $x_1+ x_2+\cdots +x_{k-1}=s$, which is $\binom{s+k-2}{k-2}$  (NB. $k\ge 2$). Thus for fixed $(i_1i_2\cdots i_k)$,
the choices of $(j_1,j_2,\cdots, j_k)$ are not more than
$$\sum_{s=1}^{[\eps n]} \binom{s+k-2}{k-2}\le \sum_{s=1}^{[\eps n]} 2^{s-2+ \eps n}< 2^{2\eps n}.$$
Here we have used the fact $\binom{m}{l}\le 2^m$ and $k\le \epsilon n$.
%
Thus
$$\#\mathcal{F}^{n,m}(\eps) \le 2^{2\eps n} \sum_{k=1}^{[\eps n]}\binom {n-1}{k-1} .$$
The last sum is bounded by $2^{n H(\epsilon)}$
where $H(\epsilon) = -\epsilon \log_2 \epsilon -(1-\epsilon)\log_2 (1-\epsilon)$, which tends to $0$ as $\epsilon \to 0$. This is given by the large deviation estimation for independent symmetric Bernoulli variables.
Thus the lemma follows.
\end{proof}

Put
$$Z(\eps) :=\bigcap_{N=1}^\infty \bigcup_{n=N}^\infty \bigcup_{m=1}^\infty \bigcup_{(\textbf{i}, \textbf{j})\in \mathcal{F}^{n,m}(\eps)} \widehat{Q}_{\textbf{i}}^{\textbf{j}}.$$
\begin{prop}\label{cor:dimzeps} $\dim_H Z(\eps)\le \eps' +\eps+ \eps \log_2 M$, where $\eps'$ is as in Lemma~\ref{lem:counting}.
\end{prop}
\begin{proof} By Lemma~\ref{lem:covering}, for each $(\textbf{i},\textbf{j})\in \mathcal{F}_k^{n,m}(\eps)$, $\widehat{Q}_{\textbf{i}}^{\textbf{j}}$ is covered by at most $M^k 2^{\eps n}$ intervals in $\mathcal{D}_{n+m}$.  As $k\le \eps n$, $M^k\le 2^{\eps \log_2 M n}$. From Lemma~\ref{lem:counting} we get
	$$
	   \sum_{n=N}^\infty \sum_{m=1}^\infty C 2^{\epsilon' n}\cdot M^k 2^{\epsilon n}\cdot  2^{-\delta(n+m)} 
	   \le C \sum_{m=1}^\infty 2^{-\delta m} \sum_{n=N}^\infty 2^{-(\delta - \epsilon' - \epsilon - \epsilon \log_2 M )n}.
	$$
	The last sum  is finite when $\delta >  \epsilon' + \epsilon + \epsilon \log_2 M $. So, the statement follows.
\end{proof}


\begin{prop}\label{prop:inzeps} For each $\eps>0$, $X_\infty\setminus Y_\infty\subset Z(\eps)$.
\end{prop}
\begin{proof} Let$\eps>0$ and let $x\in X_\infty\setminus Y_\infty$. Since $x\not\in Y_\infty$,  there is $K>0$ such that
\begin{equation}\label{eqn:slowK}
\forall n\ge 1, \quad 
\dist(T^n x, x)\ge 2^{-Kn}
\end{equation}
(such a point $x$ doesn't return back too closely to itself).
To show that $x\in Z(\eps)$, we may assume that $\eps$ is as small as we wish, since $Z(\eps)$ is increasing in $\eps$.

Let $A=1/\eps^2$. This constant $A$ will describe a tiny ball around $x$, that ``waits" for the return of $x$ infinitely many times, which are defined as follows. 
We define sequences $(i_k)_{k=1}^\infty$ and $(q_k)_{k=1}^\infty$ of positive integers as follows.
\begin{enumerate}
\item $i_1$ is the first return time after $A$ into the ball $B(x, 2^{-A})$, i.e. the minimal positive integer such that 
$$i_1\ge A,\,\, \dist(T^{i_1}x, x)< 2^{-A},$$
and let $q_1: =[-\log_2 \dist(T^{i_1}x, x)]\ge A$.
\item Once $i_k$ and $q_k\ge A$ are defined, let $i_{k+1}$ be the minimal positive integer such that 
$$i_{k+1}> i_k+ [\eps q_k], \quad \dist(T^{i_{k+1}}x, x)\le 2^{-A},$$ and
let $q_{k+1}:=[-\log_2\dist(T^{i_{k+1}}x, x)]\ge A$.
\end{enumerate}

We shall prove that there is an arbitrarily large $k$ such that
\begin{equation} \label{eqn:ilql}
i_1+\sum_{\substack{1\le l<k\\ i_{l+1}-i_l> q_l}} (i_{l+1}-i_l)< \eps i_k.
\end{equation}
This will complete the proof: putting $j_l= \min (i_{l+1}-i_l, q_l)$ for each $1\le l<k$ and $j_k=q_k$, we have
 \begin{equation}\label{X1}
 x\in Q_{i_1i_2\cdots i_k}^{j_1j_2\cdots j_k}; \quad
(i_1i_2\cdots i_k; j_1j_2\cdots j_k)\in \mathcal{F}_k(\eps).
\end{equation}
Indeed, since $j_l\le q_l$, we have
$$d(T^{i_{l}} x, x)
<2^{-q_{l}}
\le 2^{-j_l}.
$$
Thus we have checked the first point in (\ref{X1}). To check the second
point in (\ref{X1}), we are going to verify the properties (i), (ii) and (iii) for $(i_1i_2\cdots i_k; j_1j_2\cdots j_k)$:

(i) 
$j_l\le i_{l+1}-i_l$ clearly holds for each $1\le l<k$, by the definition of $j_l$; \\
\indent (ii) $k\le \eps i_k$ holds because of $i_1, i_2-i_1, \cdots, i_k-i_{k-1}\ge \eps A=1/\eps$, which implies
$$
    i_k =i_1 + (i_2 -i_1)+\cdots + (i_k - i_{k-1}) \ge \frac{k}{\epsilon};
$$ 
\indent (iii) Write
\begin{multline*}
i_k-\sum_{1\le j<k} j_l =i_1+\sum_{1\le l<k} (i_{l+1}-i_l-j_l)= i_1+\sum_{\substack{1\le l<k\\ i_{l+1}-i_l>j_l}} (i_{l+1}-i_l-j_l).
\end{multline*}
Then, by (\ref{eqn:ilql}),
$$
i_k-\sum_{1\le j<k} j_l 
\le i_1+\sum_{\substack{1\le l<k\\ i_{l+1}-i_l>j_l}} (i_{l+1}-i_l)< \eps i_k.
$$

Now let us prove that (\ref{eqn:ilql}) holds for infinitely many $k$. We start with

{\bf Claim.} For each $i_l<i\le i_l+[\eps q_l]$,
we have
\begin{equation*}
2\dist (T^i x, x)\
\ge \dist (T^{i-i_l} x, x),
\end{equation*}
and hence 
\begin{equation}\label{eqn:binding}
-\log_2 \dist (T^i x, x)\le -\log_2 \dist (T^{i-i_l} x,x)+1.
\end{equation}
Indeed, 
by (\ref{eqn:slowK}), for small $\epsilon$ we have
\begin{eqnarray*}
\dist (T^{i-i_l} x,x) &\ge& 2^{-K (i-i_l)} \ge 2^{-2K\eps q_l} > 2 \cdot 2^{-(1-\eps) q_l}
\end{eqnarray*}
As $\epsilon q_l \ge i- i_l$ and $2^{-q_l}\ge d(T^{i_l}x, x)$, it follows that
$$
\dist (T^{i-i_l} x,x) \ge 2\cdot 2^{i-i_l}\dist (T^{i_l}x, x)=2 \cdot \dist (T^i x, T^{i-i_l}x).
$$
The  claim follows from this and the triangle inequality
$$\dist (T^i x, x)\ge \dist (T^{i-i_l}x,x)-\dist (T^i x, T^{i-i_l} x).
$$

Now let
$$\xi_n=\frac{1}{n} \sum_{i=1}^n -\log_2 \dist (T^i x, x).$$
Recall that $\limsup \xi_n=+\infty$. Fix any $0<\epsilon <1/2$. There exists an arbitrarily large $N$, such that 
\begin{equation}\label{xi_N}
\xi_N> \xi_n \ (\forall 1\le n<N); \quad  
i_1\le \frac{\eps}{4}N; \quad \xi_N>8(A+2+K)/\eps.
\end{equation}
Choose $k$ such that $i_k\le N<i_{k+1}$. We shall prove that (\ref{eqn:ilql}) holds for this $k$.

For simplicity, in the sequel we 
denote 
$$ a_n(x) = - \log_2 d(T^n x, x).$$

By the claim (see (\ref{eqn:binding})), for each $l\le k$, if $i_l<t\le i_{l}+[\eps q_l]$, we have
\begin{equation}\label{E1}
\sum_{i_l\le i\le t}a_i(x) \le q_l+ \sum_{i'=1}^{t-i_l} \left(a_{i'}(x)+1\right)
\le (t-i_l)(\xi_N +1)+q_l
\end{equation}
where we have used the first inequality in (\ref{xi_N});
and if $i_l+[\eps q_l]<t< i_{l+1}$, then
\begin{equation}\label{E2}
\sum_{i_l\le i\le t}a_i(x)\le [\eps q_l](\xi_N +1)+q_l +A(t-i_l-[\eps q_l]).
\end{equation}
Here we  cut the sum into two parts, one for $i_l\le i \le i_l+[\epsilon q_l]$
for which we applied (\ref{E1}) and the other for 
$i_l+[\epsilon q_l]< i <t$ for which we used $d(T^i x, x)\ge 2^{-A}$.
Apply (\ref{E2}) to $t = i_{l+1}-1$. Since $i_{l+1}-i_l\ge [\eps q_l]> q_l/A$, we get the following consequence of (\ref{E2})
\begin{equation}\label{E3}
\sum_{i_l\le i<i_{l+1}}a_i(x)\le \left(\xi_N+1+2A\right) (i_{l+1}-i_l).
\end{equation}
Moreover, if $i_{l+1}-i_l> q_l$,
the first term on the right hand side  of  (\ref{E2}) is at most $\eps q_l (\xi_N+1)  < \eps(\xi_N+1) (i_{l+1}-i_l)$; the second term bounded by $A(i_{l+1}-i_l)$; the third term bounded by $A(i_{l+1}-i_l)$.  Thus we have
\begin{equation}\label{E4}
\sum_{i_l\le i<i_{l+1}}a_i(x) \le \left(\eps (\xi_N+1)+2A\right) (i_{l+1}-i_l).
\end{equation}

On the other hand, we have
\begin{equation}\label{E5}
\sum_{i_k\le i\le N} a_i(x) \le \left( \xi_N+A+2\right) (N-i_k) +q_k.
\end{equation}

Write
$$N\xi_N=\left( \sum_{i=1}^{i_1-1} 
+\sum_{\substack{1\le l<k \\  i_{l+1}-i_l>j_l} } \sum_{i_l\le i<i_{l+1}}
+\sum_{\substack{1\le l<k  \\ i_{l+1}-i_l\le j_l}} \sum_{i_l\le i<i_{l+1}}
+\sum_{i_k\le i\le N}\right)a_i(x).$$
From (\ref{xi_N}), (\ref{E3}), (\ref{E4}) and (\ref{E5}) we obtain
\begin{multline*}
N\xi_N\le\xi_N(i_1-1) + (\eps \xi_N+A+2) \sum_{\substack{1\le l<k\\ i_{l+1}-i_l>j_l}} (i_{l+1}-i_l)+\\
 +(\xi_N+A+2) \left(\sum_{\substack{1\le l<k\\ i_{l+1}-i_l\le j_l}} (i_{l+1}-i_l)+ N-i_k\right) + q_k,
\end{multline*}
which implies that
$$
N\xi_N\le  (\xi_N+A+2) N -(1-\eps) \xi_N \sum_{\substack{1\le l<k\\ i_{l+1}-i_l>j_l}} (i_{l+1}-i_l) + q_k.
$$
Since $q_k\le K i_k \le K N$ (see the definition of $q_k$ and (\ref{eqn:slowK})), it follows that
\begin{equation}\label{E6}
\sum_{\substack{1\le l<k\\ i_{l+1}-i_l>j_l}} (i_{l+1}-i_l)\le \frac{(A+2+K)}{(1-\epsilon)\xi_N} N <\frac{\eps}{4} N
\end{equation}
where for the last inequality we have used the third point in (\ref{xi_N}).
Since $\xi_N\ge \xi_{N-1}$ and $\xi_N>A$ (see (\ref{xi_N})), we have
$$-\log_2 \dist (T^N x,x)\ge \xi_N>A,$$
which implies
\begin{equation}\label{E7}
N\le i_k+[\eps q_k]\le i_k + K\eps i_k< 2 i_k,
\end{equation}
where the first inequality is implied by the fact that $i_{k+1}$
is the first time after $i_k+ [\epsilon q_k]$ for $x$ to go into the 
$2^{-A}$-neighbourhood of $x$. 
Finally (\ref{eqn:ilql}) follows from (\ref{E6}), (\ref{E7}) and the fact $i_1\le \eps N/4$ (see (\ref{xi_N})).
\end{proof}
\medskip

%

\begin{proof}[Proof of Theorem~\ref{thm:zerodim}]
By Proposition~\ref{prop:inzeps} and Proposition~\ref{cor:dimzeps}, we get $$\dim (X_\infty\setminus Y_\infty)=0.$$
This,  together with Proposition~\ref{lem:dimY},  implies that $\dim (X_\infty)=0$.
\end{proof}

\section{Pointwise behavior of $\prod_{k=0}^{n-1} |\cos \pi 2^k x|$}
Let us consider the particular case $f_{0}(x) = \log |\cos \pi x|$
and look at pointwise behavior of its Birkhoff averages and compute
$\alpha (0)$, $\beta(0)$.

\begin{theorem}\label{thm:MAf0} We have $\alpha(0) = - \log 2$ and $\beta(0) =0$. There are only
three possible limits for
\begin{equation}\label{lim-L1/2}
     \lim_{n\to\infty}\frac{1}{n}  \sum_{k=0}^{n-1} \log |\cos \pi 2^k x|.
\end{equation}
The limit 0 is attained at $x=0$, the limit $-\infty$ is attained at dyadic rational points $\frac{2k+1}{2^m}$,  and $-\log 2$
is the other possible limit which is attained by a set of full Lebesgue
measure.
	\end{theorem}

\begin{proof}
Our discussion is based on the following relation
\begin{equation}\label{Eq_1/2}
\frac{1}{n} \sum_{k=0}^{n-1}  \log  |\cos \pi 2^k x| =   -\log 2 +  \frac{\log |\sin \pi 2^n x|}{n} - \frac{\log |\sin \pi x|}{n}
\end{equation}
(for $x\not=0$), which is the consequence of the fact that $\log |\cos \pi|$ is a coboundary, i.e.
$$
\log |\cos \pi x| = \log |\sin \pi 2x| - \log |\sin \pi x| -\log 2.
$$
Let $\mu$ be a $\tau$-periodic invariant measure with $\tau \ge 2$. The support of the measure contains neither $0$, the singularity of $-\log |\sin \pi x|$, nor $1/2$,
the singularity of $-\log |\cos \pi x|$. So,
by  the Birkhoff theorem (or a direct argument), letting $n\to \infty$ in
(\ref{Eq_1/2})  we get 
$$
\int \log |\cos \pi x| d\mu(x)= \frac{1}{\tau}  \sum_{k=0}^{\tau-1} \log  |\cos \pi (2^k x_0)|  =  - \log 2
$$
where $x_0$ is a point in the $\tau$-cycle supporting $\mu$.
Recall that $0$ is the unique fixed point and  $\int \log |\cos \pi x| d\delta_0(x) =0$. Thus we have proved
$$
\alpha(0)  =  \inf_{\mu \ \mbox{\rm \small periodic}} \int  \log |\cos \pi x| d\mu = -\log 2
$$
and
$$
\beta(0) = \int \log |\cos \pi x| d\delta_0(x) =0.
$$

Now consider the possible value of the limit (\ref{lim-L1/2}). The limit is equal to $0$ when  $x=0$, and it is equal to $-\infty$ if $x=\frac{2k+1}{2^m}$ because $2^{n} x = 0 \mod 1$ for all $n> m$
so that the RHS remains $-\infty$ for all $n>m$. 
Suppose now that $x$ is not dyadic rational, or equivalently  $2^n x \not =0$ for all $n\ge 0$ and that the limit (\ref{lim-L1/2}) exits. We claim that the limit is equal to $-\log 2$, or equivalently, 
by (\ref{Eq_1/2}),
$$
\lim_{n \to \infty}   \frac{\log |\sin \pi 2^n x|}{n}=0.
$$ 
Otherwise, let $-a<0$ be the limit. Then $|\sin \pi 2^n x| \approx e^{- an }$.
But $2^n x = \sum_{j=1}^\infty \frac{x_{n+j-1}}{2^j}$. Thus for any  $n$
large enough, the word $x_{n}x_{n+1}\cdots$ has a long prefix of $0^{bn}$.
It follows that $2^n x=0$ when $n$ is large enough, a contradiction to   $|\sin \pi 2^n x| \approx e^{- an }$.
	\end{proof}

We have actually proved that if we take limsup in (\ref{lim-L1/2}),
we get still only three possible values. However, if we consider liminf, any value
$a <-\log 2$ is possible.


Let us compute the pressure function $p_0(\cdot)$. First we have obviously  
$\int f_0 d\delta_0 =0$ It follows that $p_0(t) \ge 0$ for all $t$. 
Secondly, observe that any invariant measure $\mu$ can be written as $a \delta_0 + (1-a)\mu'$
where $0\le a\le 1$ and $\mu'$ is an invariant measure which does not have $0$ as atom. The measure $\mu'$ does not have dyadic rational points $\frac{k+1}{2^{m}}$ as atoms either. So, by Theorem  \ref{thm:MAf0} and the ergodic theorem, we obtain 
   $$
      \int f_0 d\mu = - (1-a) \log 2. 
   $$ 
   It follows that
$$
       \forall t \in \mathbb{R}, \ \ \ p_0(t) = \max \{(1- t)\log 2,0\}.
$$ 
Then
$$
     p_0^*(\alpha) = |\alpha|
     \ \  {\rm for } \   \alpha \in [-\log 2, 0];
     \quad p_0^*(\alpha) = -\infty\  {\rm for \ other } \ \  \alpha.
$$

\section{Multifractal analysis of $\prod_{k=0}^{n-1}|\cos \pi (2^k x +c)|$}
In this last section, 
we are going to prove  Theorem B announced  in Introduction.

The difficulty of proof is due to the fact that $\log |\cos \pi (2^j x +c)|$
is atypically large (in absolute value) when $2^j x$ approaches the singularity $b:=1/2-c$.

By the way, we point out that the set $E_c(\alpha)$ is defined by the existence of the limit and it can be empty.
But if we consider  the set $$
  E_c^-(\alpha):=\left\{x \in \mathbb{R}/\mathbb{Z}: \liminf_{n\to\infty}\frac1n \sum_{i=0}^{n-1}f_c(T^ix)\le\alpha\right\}
  \quad (\alpha \le 0),
$$ 
we have  $\dim_H E_c^-(\alpha) >0$ for all $\alpha$. This is a contrast to the point
(4)  of Theorem B when $\alpha(c)>-\infty$ and  $\alpha <\alpha(c)$ for which $E_c(\alpha)=\emptyset$. To see this just consider sequences with infinitely often repeating patterns: $b_1\cdots b_{[|\alpha| n]}$ where $(b_k)$ is the dyadic expansion of the singular point $b$. 

\medskip

Here is the idea for proving Theorem B. 
We shall decompose $E_c(\alpha)$ into three parts: the set of points whose orbit  
does not come  close to the singularity $b$, the set  of points whose orbit comes close to $b$ with infinitely many good times and the set of points whose orbit comes  close to $b$ with only finitely many good times. The idea of considering good times comes from \cite{PRL2007}. The first part, called hyperbolic part, can be treated by the classical thermodynamical formalism and the key is to cover the hyperbolic part by Markov subsystems. As we shall see, the third part is a small set in the sense of Hausdorff dimension.   


The three parts are separately dealt with in the following subsections. They are followed by a synthesis of proof. 




\subsection{Markov subsystems}
Let us reserve the name {\em Markov subsystem} for a continuous map $F: \bigcup_{j=1}^m J_j\to \T$ with the following properties:
\begin{itemize}
	\item $J_1,J_2, \cdots, J_m$ are closed subintervals of $\T\setminus \{b\}$ with pairwise disjoint interiors;
	\item for each $1\le j\le m$, there is a positive integer $n_j$ such that $F|J_j:=T^{n_j}|J_j$ is a homeomorphism from $J_j$ onto a subinterval
	$T^{n_j}(J_j)$ of $\T$ and $F^i(J_j)\not\ni b$ for each $0\le i<n_j$;
	\item for each couple $1\le j,k\le m$, either $F(J_j^o)\cap J_k^o=\emptyset$ or $F(J_j)\supset J_k$.
\end{itemize}
As usual, a Markov subsystem is said to be  {\em irreducible} if for all $1\le j, k\le m$, there exists $j_0=j,j_1,\cdots, j_s=k$ such that $F(J_{j_i})\supset J_{j_{i+1}}$ for each $0\le i<s$; and an irreducible Markov subsystem is said to be {\em aperiodic} if there is a positive integer $N$ such that for any $1\le j, k\le m$, there exists $j_0=j, j_1, \cdots, j_N=k$ such that $F(J_{j_i})\supset J_{i+1}$ for each $0\le i<N$. 
The integers $n_j$ in the above definition are called {\em induced times}. 

A Markov subsystem is irreducible if and only if its maximal invariant set $\bigcap_{n=1}^\infty F^{-n}(J_1 \cup \cdots \cup J_m)$ is $F$-transitive. An irreducible Markov subsystem is aperiodic if there is $j$ such that $F(J_j)\supset J_j$.

Let $T:\bigcup_{j=1}^m J_m\to \T$ be a Markov subsystem with a maximal invariant set $K$. Notice that here the induced times are all assumed to be $1$.  Let $\mathcal{M}_K$ denote the collection of all invariant probability Borel measures supported in $K$. We define the pressure 
of $f_c$ associated to the subsystem $T: K \to K$ by 
$$p(t;K):=\sup\left\{h_\mu+t\int f_c d\mu: \mu\in \mathcal{M}_K\right\},$$
and its Legendre transform $p^*(\alpha;K)=\inf_{t\in\R} (p(t;K)-t\alpha)$. 
Also define $$
\alpha_K:=\inf_{\mathcal{M}_K} \int f_c d\mu, \qquad \beta_K:= \sup_{\mathcal{M}_K} \int f_c d\mu.$$

Since $f_c$ restricted on $K$ is Lipschitz continuous, the set
$K \cap E(\alpha)$ supports a Gibbs measure.  The following is well-known, see \cite{BPS1997,F1994}.

\begin{lemma}\label{lem:clatheomo} Suppose that $T:\bigcup_{j=1}^m J_m\to \T$  is an irreducible aperiodic Markov system with a maximal invariant set $K$. 
	Then 
	\begin{itemize}
		\item[{\rm (i)}] 	for any  $\alpha \in (\alpha_K, \beta_K)$,  we have $$\dim (K\cap E(\alpha))=\frac{p^*(\alpha;K)}{\log 2};$$
	    \item[{\rm (ii)}] for $\alpha>\beta_K$ and for $\alpha<\alpha_K$, $K\cap E(\alpha)=\emptyset$;
	    \item [{\rm (iii)}] for $\alpha=\alpha_K$ and for $\alpha=\beta_K$,
	    $\dim (K\cap E(\alpha))\le \frac{p^*(\alpha;K)}{\log 2}.$  
	\end{itemize}

\end{lemma}

Our main technical result is the following proposition. 
For each $\delta>0$, let
$$K_\delta:=\{x\in \T: d(T^n(x), b)\ge \delta, \forall n\ge 0\}.$$
It is a closed $T$-invariant set. 

\begin{prop}\label{prop:mixing}
	Assume $0<c<1$. Then for each $\delta>0$ there is an irreducible aperiodic Markov subsystem $T:\bigcup_{j} J_j\to \T$ whose maximal invariant set contains $K_\delta$, and doesn't contain $b$.
\end{prop}

The proof of the proposition relies on the following lemma. Let us first introduce a notation. 
Fix $x \in \mathbb{T}$. Its forward orbit and finite forward orbits are denoted by
$$
O^+(x) :=\{x, Tx, T^2x, \cdots\}, \quad O^+_n(x) :=\{x, Tx, \cdots, T^{n-1} x\} \ (\forall n \ge 1);
$$
and its backward orbit and finite backward orbits are denoted by
$$
O^-(x) :=\bigcup_{k=1}^\infty T^{-k}x,  \quad O^-_n(x) :=\bigcup_{k=1}^n T^{-k}x \ (\forall n \ge 1).
$$
For example, for $x=\frac{1}{2}$,  we have 
$$
O^+_n(\frac{1}{2})  = O^+(\frac{1}{2}) = \{0, 1/2\} \ \ (\forall n \ge 2);
$$
$$
 T^{-n}\frac{1}{2} = \left\{\frac{2k+1}{2^{n+1}}: 0\le k<2^n\right\}\ (\forall n \ge 1).
$$
For $x=0$ we have $O^+(0)=O^+_n(0) = \{0\}$ ($\forall n\ge 1$) and 
$$
O^{-}(0) = \left\{0, \frac{1}{2}\right\} \cup \bigcup_{n=1}^\infty \left\{\frac{2k+1}{2^{n+1}}: 0\le k<2^n\right\}.
$$

In the remaining of this subsection, we shall prove Proposition~\ref{prop:mixing}. 
Notice that
the following lemma  fails in the case $c=0$, because the conclusion fails for $x=0$. But the conclusion does hold for all $x \in \mathbb{T}\setminus \{0\}$ when $c=0$, as the proof shows.
The lemma states that any point $x$ admits a preimage $y$ such that its backward orbit and one of its finite forward orbits never encounter the
singularity.

\begin{lemma}\label{lem:preimage}
	Suppose $0<c<1$ (recall that $b:=\frac{1}{2}-c \not= \frac{1}{2}$).  For any  $x\in \T$, there exist an integer $n:=n(x)\ge 1$ and a point $y\in T^{-n} x$ such that $b \not\in O^+_n(y) \cup O^-(y)$.  
\end{lemma}
\begin{proof} Let $X$ denote the set of all $x\in \T$ with the following property: 
	$$
	\forall n\ge 1,  \forall y\in T^{-n}x, \exists 0\le j <n\ {\rm s.t.} \ T^jy=b \ {\rm or}\  \exists m\ge 1 \ {\rm s.t.}\ T^mb=y.
	$$ 
	Then $x \in X^c$ means
	$$
	\exists n\ge 1,  \exists y\in T^{-n}x, \forall 0\le j <n, \ T^jy\not=b \ {\rm and}\  \forall m\ge 1, T^mb\not=y.
	$$ 
	So, what we have to prove is $X=\emptyset$. First note that
	\begin{equation}\label{eqn:Xbackinv}
	T^{-1}(X)\setminus \{b\}\subset X.
	\end{equation}
	
	{\em First step: each $x\in X$ is periodic}. We distinguish two cases. Case 1: $x\not=T(b)$. Consider the two preimages $y_1, y_2$ of $x$ under $T$,
	which are  not necessarily equal to $b$. There exist two integers  $m_1, m_2\ge 0$ such that $y_i=T^{m_i}(b)$, $i=1,2$. Then $m_1\not=m_2$ and $x=T^{m_1+1}(b)=T^{m_2+1}(b)$ is periodic of period $|m_2-m_1|$. Case 2: $x=T(b)$. Consider  $y:= b+\frac{1}{2}$, the preimage of $x$ under $T$ other than $b$. Then $y\in X$, by (\ref{eqn:Xbackinv}). If $y\not=T(b)$ then  $y$ is periodic by Case 1 and hence so is $x$. If  $y=T(b)$ then $x=T(b)=T^2(b)$ is fixed.
	
	{\em Second step:  $T(b)$ is the only possible member of $X$}. To see this, let $x\in X$ and 
	let $y$ be the preimage of $x$ under $T$ which is not in the orbit of $x$. Then $y$ is not periodic, hence $y\not\in X$. Therefore, $y=b$ and $x=T(b)$. 
	
	{\em Final step: $X=\emptyset$}. Arguing by contradiction, assume $X\not=\emptyset$. Then by the preceding two steps, we have $X=\{0\}$ and $T(b)=0$. As $0<c<1$, this forces $b=0$. However, $1/2\in T^{-1}(0)$ and $T^m(0) \not= 1/2$ for any $m\ge 0$. This means $0\not\in X$, a contradiction!
\end{proof}

\begin{proof}[Proof of Proposition~\ref{prop:mixing}]
	Fix $\delta>0$ small enough so that $K_\delta$ contains two periodic points
	with coprime periods. 
	
	It suffices to show that $K_\delta$ is contained in a compact $T$-invariant set $E\subset \T\setminus \{b\}$ such that $T|_E$ is transitive. Indeed, let $N_0$ be a large integer and let $\mathcal{J}$ be the union of all closed dyadic intervals of length $2^{-N_0}$ which intersect $E$ ($N_0$ large enough so that these dyadic intervals do not contain $b$). Then consider the Markov subsystem $T|_{\mathcal{J}}:\mathcal{J}\to \T$. 
	The transitivity of $T|_E$ implies that $T|_{\mathcal{J}}$ is irreducible, and the fact that  $K_\delta$ as well as $E$ contains 
	 two periodic points
	with coprime periods implies that this Markov subsystem $T|_{\mathcal{J}}$ is aperiodic.
	
	To this end, we start with choosing an arbitrary periodic point
	$v$ of $T$
	outside $K_\delta$ such that $b$ is not on its forward orbit, i.e. $v\not\in K_\delta$ with $b \not\in O^+(v)$. Let $s$ be the period of $v$. Take a small open dyadic interval $V\ni v$ such that $b \not\in T^j(\overline{V})$ for all $0\le j<s$. For every $k\ge 1$, let $V_k$ denote the component of $T^{-ks} (V)$ which contains $v$.
	
	Applying Lemma~\ref{lem:preimage}, for any $x\in K_\delta$ we obtain a time $n \ge 1$ and a point $y\in T^{-n} x$ such that $b \not\in O_n^{+}(y) \cup O^{-}(y)$. Pull back $y$ further to meet $V$, say $\ell$ times. We can obtain  $n(x):=n+\ell\ge 1$ and $\rho(x)>0$  such that $T^{-n(x)}(B(x,\rho(x)))$ has a component $U(x)$ in $V$ with $T^j(\overline{U(x)})\not\ni b$ for each $0\le j<n(x)$. The existence of $\rho(x)$ is ensured by the fact $b \not\in O_n^{+}(y) \cup O^{-}(y)$.
	By compactness of $K_\delta$, there is an integer $N$ such that for each closed dyadic interval $I$ of length $2^{-N}$ which intersects $K_\delta$, there is $n'_I\ge 1$ and $U'_I\subset V$ such that $$
	T^j(U'_I)\not\ni b\ (\forall 0\le j<n'_I);  \quad T^{n'_I}(U'_I)=I.$$ Let $\mathcal{I}$ denote the collection of all the dyadic intervals of length $2^{-N}$ which intersect $K_\delta$.
	
	Again, as above,  apply Lemma~\ref{lem:preimage} to the couple $(v, I)$ instead of the couple $(x, V)$. For each $I\in \mathcal{I}$, there exist $m_I \in  \mathbb{N}$ and $y_I\in I^o$ such that $T^{m_I}(y_I)=v$ and $T^j(y_I)\not=b$ for each $0\le j<m_I$. So there exists a (large) positive integer $k$ such that for each $I$, there exists $W_I\subset I$ such that 
	$$
	     T^{j}(W_I)\not\ni b \ (\forall 0\le j<m_I); \quad T^{m_I}(W_I)=V_k.
	$$ 
	Fix the choice of $k$ now. As $K_\delta$ is $T$-invariant and $v\not\in K_\delta$, choosing $k$ large enough, we can ensure that $W_I\cap K_\delta=\emptyset$. Now take $M$ large so that $d(W_I, K_\delta)>2^{-M}$ for each $I$, and let $J_1,J_2,\cdots, J_\nu$ denote all the dyadic intervals of length $2^{-M}$ which intersect $K_\delta$.
	
	In the next step, for each $I\in \mathcal{I}$, we choose a dyadic interval $U_I$ and $n_I\equiv n'_I\mod s$, so that $U_I\subset V_k$, $T^{n_I}(U_I)=I$ and $T^j(U_I)\not\ni b$ for each $0\le j<n_I$. We can arrange that these intervals $U_I$ are pairwise disjoint, since for each $I$, we can choose $U_I$ arbitrarily closed to $v$.
	
	Finally, define a Markov subsystem $F$ as follows:\\
	\indent \ \ \ $F|_{J_j}=T^{M-N}|_{J_j}$ for each $j$, \\
	\indent \ \ \ $F|_{W_I}=T^{m_I}|_{W_I}$, \ $F|_{U_I}=T^{n_I}|_{U_I}$ for each $I\in\mathcal{I}$. \\
	Then $F$ is irreducible and hence its maximal invariant $E'$ is $F$-transitive. Define $$
	E=\bigcup_{n=0}^\infty T^n(E').$$
	 This is a $T$-invariant compact subset disjoint from $b$ and containing $K_\delta$.
\end{proof}

\subsection{Proof of $\dim E_c^{hyp}(\alpha) = \frac{p_c^*(\alpha)}{\log 2}$}
For the subsystem $T: K_\delta \to K_\delta$, we define its pressure 
function by
$$
p(t; \delta):=\sup\left\{ h_\mu+ t\int f_c d\mu: \ \mu\in\mathcal{M}_T(\delta)\right\},$$
where
$$ \mathcal{M}_T(\delta)=\left\{\mu\in \mathcal{M}_T: \supp(\mu)\subset K_\delta\right\}.$$
We also define the extremal values
$$\alpha_\delta(c):=\inf_{\mu\in \mathcal{M}_T(\delta)} \int_{\mathbb{T}} f_c(x) d\mu, \quad \beta_\delta(c):=\sup_{\mu\in\mathcal{M}_T(\delta)} \int_{\mathbb{T}} f_c(x) d\mu.$$
Since both $\alpha(c)$ and $\beta(c)$ can be defined by periodic measures
(cf. Proposition \ref{prop:ab}), 
we have
\begin{equation}
     \alpha_\delta(c) \downarrow \alpha(c), \quad \beta_\delta(c) \uparrow \beta(c) \ \ \ {\rm as} \ \delta \downarrow 0. 
\end{equation}

\begin{prop}\label{lem:ptinf} For any $\alpha\in (\alpha(c),\beta(c))$,
	\begin{equation}\label{eq:ptinf}
	\log 2\cdot \dim E_c^{hyp}(\alpha) = \inf_{t\in \R} (p(t)-t\alpha)=\lim_{\delta\searrow 0} \inf_{t\in \R} (p(t; \delta)-t\alpha),
	\end{equation}
\end{prop}
\begin{proof} Let I,II and III denote the three terms appearing in (\ref{eq:ptinf}). 
	
	 We first prove the equality $\textrm{II}=\textrm{III}$. Since $p(t)\ge p(t; \delta)$ for each $\delta>0$ and $t\in \R$, $\textrm{II}\ge \textrm{III}$. To prove the opposite inequality, assume by contradiction that $\textrm{II}>\textrm{III}$. 
	Then there exist $\eta>0$, $\delta_n>0$ and $t_n\in\R$ ($n=1,2,\ldots$) with $\delta_n\to 0$, such that 
	\begin{equation}\label{eqn:pdeltantn}
	p(t_n;\delta_n)-t_n\alpha+\eta< p(t)-t\alpha, \text{ for all } t\in \R.
	\end{equation}
Let us show that $\{t_n\}_{n=1}^\infty$ is bounded. 
Let $\eps>0$ be such that $[\alpha-\eps, \alpha+\eps]\subset (\alpha(c), \beta(c))$. Then there exists $\delta_*>0$ such that for any $\delta\in (0,\delta_*)$,
$[\alpha-\eps, \alpha+\eps]\subset (\alpha_\delta(c),\beta_{\delta}(c))$. 
For each $\mu\in\mathcal{M}_T(\delta)$, we have
$$p(t; \delta)-t\alpha\ge  t\left(\int f_c d\mu-\alpha\right).$$
Thus 
\begin{equation}\label{eqn:pdeltatlow}
p(t; \delta)-t\alpha\ge \left\{
\begin{array}{ll}
|t|(\alpha-\alpha_\delta(c)) &\mbox{ if } t< 0,\\
t(\beta_\delta(c)-\alpha) &\mbox{ if } t\ge 0.
\end{array}
\right.
\end{equation}
It follows that for each $\delta\in (0,\delta_*)$, 
$$p(t; \delta)-t\alpha\ge |t|\eps.$$
Then, by (\ref{eqn:pdeltantn}), we obtain that $t_n$ is bounded. 
Passing to a subsequence, we may assume that $t_n$ converges to some $t_*\in\R$ and that $t_n$ is monotone in $n$. 
Suppose first that $t_n$ is monotone increasing, so that $t_n\le t_*$ for all $n$.  
%
Since $p(t; \delta)$ is monotone decreasing in $t$ and $p(t_*;\delta_n)\to p(t_*)$, we obtain 
$$p(t_n; \delta_n)\ge p(t_*; \delta_n)\to p(t_*),$$
which contradicts (\ref{eqn:pdeltantn}) for large $n$ with $t=t_*$.
Now suppose that $t_n$ is decreasing.
Since $t\mapsto p(t; \delta)$ is continuous in $t$ and $\delta\mapsto p(t; \delta)$ is monotone decreasing in $\delta$, we have 
$$\limsup_{n\to\infty} p(t_n;\delta_n)\ge \limsup_{n\to\infty} p(t_n; \delta)=p(t_*;\delta),$$
for any $\delta>0$. Therefore $$\limsup_n (p(t_n;\delta_n)-t_n\alpha)\ge \sup_{\delta>0} p(t_*;\delta)-t_*\alpha=p(t_*)-t_*\alpha,$$
contradicting (\ref{eqn:pdeltantn}) again. 
This proves that $\textrm{II}\le \textrm{III}$, and hence the equality $\textrm{II}=\textrm{III}$.   
  
Now we prove the equality $\textrm{I}=\textrm{III}$.
	For each $n\ge 1$, by Proposition \ref{prop:mixing}, there is an irreducible aperiodic Markov subsystem $T: \mathcal{J}_n\to \T$ with maximal invariant set $E_n\supset K_{1/n}$ and $E_n\not\ni b$. Let $\tilde{p}_n(t)$ denote the pressure of $(T|_{E_n}, f_c)$. Then  obviously 
	$$
	     p(t;1/n)\le \tilde{p}_n(t)\le  p(t),
	$$ 
	So, 
	\begin{equation}\label{p*}
	\textrm{III}=\lim_{n\to\infty} \inf_{t\in \R} (\tilde{p}_n(t)-t\alpha)
	\end{equation}
	for each $\alpha\in (\alpha(c),\beta(c))$. Notice that
	$$\inf_{\mathcal{M}_{E_n}} \int f_c d\mu \to \alpha(c);\qquad \sup_{\mathcal{M}_{E_n}} \int f_c d\mu\to\beta(c).$$
	Thus, by Lemma~\ref{lem:clatheomo}, for $n$ sufficiently large, we have
		\begin{equation}\label{h*}
		\log 2\cdot \dim (E(\alpha)\cap E_n) =\inf_{t\in\R} (\tilde{p}_n(t)-t\alpha).
			\end{equation}
	As $\dim(E^\hyp(\alpha))=\lim_n \dim (E(\alpha)\cap E_n)$, the equality $\textrm{I}=\textrm{III}$ follows from (\ref{p*}) and (\ref{h*}).
\end{proof}

\subsection{Good times and Bad points}
The idea of using the notion of good time  is inspired by \cite{PRL2007}. Using this notion, we distinguish "good points" from "bad points". We shall
show that the set of bad points is small. 

Given $\frac{1}{2}>\delta>\delta'>0$, a positive integer $n$ is called a {\em $(\delta,\delta')$-good time} of $x\in \T$, if the following hold:
\begin{itemize}
	\item $\dist(T^n(x),b)<\delta'$;
	\item  $b\not\in T^j(J)$ for each $0\le j<n$ where $J$ is the component of $T^{-n}(B(b,\delta))$ containing $x$.
\end{itemize}
Let $$\mathcal{B}(\delta,\delta')=\{x\in \T: \omega(x)\ni b,\mbox{ but } x \mbox{ has no } (\delta,\delta')\!-\!\mbox{good time}\}.$$

\begin{prop}\label{lem:dimB} For any $\eps>0$ and any $C>1$, there exist $0<\delta'<C\delta'<\delta<1/2$ such that
	$$\dim (\mathcal{B}(\delta, \delta')) <\eps.$$
\end{prop}
\begin{proof} We distinguish two cases: $b$ is a periodic point or not.
	
	{\bf Case 1.} {\em $b$ is a periodic point of $T$}. Let $s$ denote the period of $b$. Choose a small constant $\delta>0$ so that
	\begin{equation}\label{eqn:bnotret}
	\forall j\ge 0, \ \ T^j(b)\not\in B(b, \delta) \mbox{ if } s\not| j.
	\end{equation}
	Take an arbitrary $\delta'\in (0,\delta)$. It is clear that $b\in \mathcal{B}(\delta, \delta')$. We shall prove that
	$\mathcal{B}(\delta, \delta')=\{b\}$ by showing that any $x\in \T\setminus \{b\}$ with $\omega(x)\ni b$ has a $(\delta,\delta')$-good time. 
	
	
	To this end, let $n=n(x)\ge 1$ denote the minimal positive integer such that
	$$\dist(T^n(x),b)< \min (\dist(x,b), \delta').$$ Let us prove that $n(x)$ is a $(\delta,\delta')$-good time of $x$. Arguing by contradiction, assume that this is not the case.  Then there exists $0\le k<n$ such that $T^k(J)\ni b$, where $J$ denotes the component of $T^{-n}(B(b,\delta))$ containing $x$. Hence $$
	T^{n-k}(b)\in T^{n-k}(T^kJ)=B(b,\delta).$$
	By (\ref{eqn:bnotret}), we must have $n-k=ms$ for some  integer $m\ge 1$.
	
	Notice that for any $0\le i \le n$,
	$T^i J$ is a ball of radius $2^{-n+i}\delta$ in $\mathbb{T}$ centered at some point belonging to $T^{-(n-i)}(b)$.  
	Since $b \in T^k(J)$ and $k=n -ms$, we have $T^k(J) =B(b,2^{-ms}\delta)$.
	
	Also notice that $T^{j-i}: T^i J \to T^j J$ is bijective for any $0\le i <j \le n$.  Since both $b$ and $T^kx$ belong to $T^k J= T^{n-ms}J$, 
	$T^{ms}$ maps $[b, T^k(x)]$ onto $[b, T^n(x)]$ bijectively, so we have $$\dist(T^k x, b)=2^{-ms} \dist (T^n x, b)< \dist (T^n x,b).$$
	By the minimality of $n$, this implies that $k=0$. It follows that
	$n=ms$, $B(b, 2^{-ms} \delta) = J \ni x$ and 
	$$  \forall \ 0\le i<m, \ \ \dist(T^{(i+1)s}(x), b)=2^s \dist(T^{is} (x), b).$$
	In particular, $\dist (T^n x,b)=2^{ms} \dist(x, b)>\dist (x,b)$, a contradiction!

	{\bf Case 2.} {\em $b$ is not a periodic point of $T$}. Let $K$ be a large positive integer to be determined, let $\delta>0$ be small enough so that
	\begin{equation}\label{choice-delta}
	\dist(x, b)<\delta\Rightarrow \dist(T^j(x), b)\ge \delta, \forall 1\le j<K.
	\end{equation}
	Such a $\delta>0$ does exist for any given $K$, because $b$ is not periodic. Otherwise, for any $m\ge 1$ there exist $x_m \in B(b, 1/m)$
	and $1\le j_m <K$ with $d(T^{j_m}x_m, b)<1/m$. Assume $j_m =\tau$
	for infinitely many $m$. Then $T^\tau b =b$, a contradiction.
	
	Let $\delta':=2^{-K+1}\delta$.
	For each $x\in \mathcal{B}(\delta,\delta')$ fixed, let 
	$$
	\mathcal{N}_x:=\{n\ge 1: \dist(T^n(x), b)<\delta'\}.
	$$
	This set is infinite because $b\in \omega(x)$.
	Then define a map
	$\mathcal{A}_x:\mathcal{N}_x\to \mathcal{N}_x\cup\{0\}$ as follows: for $n\in\mathcal{N}_x$, define
	$$\mathcal{A}_x(n):=\max\{0\le j<n: J_j\ni b\}$$
	where
	$J_j$ denotes the component of $T^{-(n-j)}(B(b,\delta))$ containing $T^j(x)$, $0\le j< n$.
	Such an integer $\mathcal{A}_x(n)$ exists for otherwise, $n$ would be a $(\delta,\delta')$-good time of $x$. Note that $n-\mathcal{A}_x(n)\ge K$
	because of our above choice of $\delta$ (see (\ref{choice-delta})), so that
	$$|J_{\mathcal{A}_x(n)}| = 2^{\mathcal{A}_x(n)} \cdot 2^{-n}\cdot 2 \delta \le 2^{-K}\cdot 2\delta=\delta'.
	$$
	So,  $\mathcal{A}_x(n)\in\mathcal{N}_x\cup\{0\}$ and the map $\mathcal{A}_x$ is well defined.
	
	So, to each $x\in\mathcal{B}(\delta,\delta')$ and to each $n\in \mathcal{N}_x$ is associated a sequence $n=n_0>n_1>\cdots> n_\ell=0$ with $n_j-n_{j+1}\ge K$ for each $0\le j<\ell$, which is defined by  $n_{j+1}=\mathcal{A}_x^j(n)$.
	On the other hand, given a finite sequence $\textbf{n}: n_0>n_1>\cdots> n_\ell=0$ with $n_j-n_{j+1}\ge K$, the set
	$$J_{\textbf{n}}:=\{x\in \mathcal{B}(\delta,\delta'): n_0\in\mathcal{N}_x, n_{j+1}=\mathcal{A}_x(n_j) \mbox{ for each }0\le j<\ell\}$$
	is either empty, or contained in a component of $T^{-n_0}(B(b,\delta))$.  For each $N$, the collection of all $J_{\textbf{n}}$ with $n_0\ge N$ form a $2^{-N}$-covering of the set $\mathcal{B}(\delta,\delta')$.
	For each positive integer $m$, the number of $\textbf{n}$ with $n_0=m$ is bounded from above by
	$$\sum_{1\le \ell \le [m/K]} \binom{m-1}{\ell-1} =O( e^{\eps m}),$$
	provided that $K$ was chosen large enough. Thus for each $\alpha>\eps$,
	$$\sum_{\textbf{n}: n_0\ge N} \diam (J_{\textbf{n}})^\alpha\le \sum_{m=N}^\infty 2^{-m\alpha }\#\{\textbf{n}: n_0=m\}=O(1).$$
	Therefore $\dim (\mathcal{B}(\delta, \delta'))\le \eps.$
\end{proof}

\subsection{Points passing by the singularity with good times}

Let $\alpha \in \mathbb{R}$. For any $\delta>0$ and $\delta'>0$, define
	$$E_*(\alpha;\delta,\delta'):=
	\{x\in E(\alpha):
 x \mbox{ has infinitely }
\mbox{ many }(\delta,\delta')-\mbox{good times}
\}$$
\begin{prop}\label{cor:dimE*} For any $\eps>0$, there exists $C>1$ 
such that for any $0<\delta'<C\delta'<\delta<1/2$, any $\alpha, t\in\R$, we have 
	\begin{equation}\label{eq:E*}
	\dim (E_*(\alpha;\delta,\delta'))\le \frac{1}{\log 2}(p(t)-t\alpha+\eps).
	\end{equation}
	It is understood that $E_*(\alpha;\delta,\delta')=\emptyset$ 
	if the right hand is negative.	
\end{prop}

To prove the proposition, we need the following lemma.
For any $\delta>0$, let $Z_{n}(\delta)$ denote the set of points $z$ with the following property: 
\begin{itemize}
\item $\dist(z,b)<\delta'$ and $T^n(z)=b$;
\item 
$b\not \in T^j(J)$ for all $0\le j<n$, where $J$ denotes the component of $T^{-n}(B(b, \delta))$ which contains $z$.
\end{itemize}
So, for each $z\in Z_n(\delta)$, $n$ is a $(\delta,\delta')$-good time of $z$ for any $0<\delta'<\delta$. 
Let
$$\Lambda(t;\delta):=\limsup_{n\to\infty}\frac{1}{n}\log \sum_{z\in Z_{n}(\delta;\delta')} e^{tS_n f_c(z)}.$$
\begin{lemma} \label{lem:Lambdat}
	For each $t\in\R$ and any $\eps>0$, there exists $C>1$ such that if $0<\delta'<C\delta'< \delta<1/2$, then $$\Lambda(t;\delta)\le p(t)+\eps.$$
\end{lemma}

\begin{proof}
	Fix $\eps>0$ and let $C>1$ be large enough such that
	\begin{equation}
	d(y,y')\le  \frac{2}{C}\dist(y, b)\Rightarrow |f_c(y')-f_c(y)|<\frac{\eps}{|t|+1}.
	\end{equation}
	%
	Let $n$ be large such that $2^{-n+1}\delta<\delta'$. For each $z\in Z_{n}(\delta,\delta')$, let $J_z$ denote the component of $T^{-n}(I)$ which contains $z$, where $I=\overline{B(b,\delta)}$. Then $J_z\subset \overline{B(b,2\delta')}$.
	Let $D=\bigcup_{z\in Z_n(\delta,\delta')} J_z$ and let $F=T^n|D:D\to I$ is a Markov map. Since $\bigcup_{j=0}^{n-1} T^j(D)$ is a compact set which does not contain $b$, we have
	$$p(t)\ge \lim_{k\to\infty} \frac{1}{kn} \log \sum_{y\in F^{-k}(b)} e^{tS_{kn} f_c(y)}.$$

	For any $z\in Z_{n}(\delta,\delta')$, $x\in J_z$ with $T^n(x)\in B(b,2\delta')$ we have
	$$\frac{\dist(T^jx, T^jz)}{\dist(T^jz, b)}\le \frac{\dist(T^jx, T^j z)}{d(T^jz,\partial T^j J_z)}=\frac{2\dist(T^n x, b)}{|I|}< \frac{2\delta'}{\delta}\le \frac{2}{C},$$
	hence $$|tS_n f_c(x)-tS_n f_c(z)|<n\eps.$$
	Therefore
	$$\sum_{y\in F^{-k-1}(b)} e^{tS_{(k+1)n}f_c(y)} \ge e^{-n\eps} e^{n\Lambda_n} \sum_{y\in F^{-k}(b)} e^{tS_{kn}f_c(y)},$$
	where
	\begin{equation}\label{eqn:dfnn}
	\Lambda_n:=\frac{1}{n}\log \sum_{z\in Z_{n}(\delta,\delta')} e^{tS_n f_c(z)}.
	\end{equation}
	Consequently, $p(t)\ge \Lambda_n-\eps.$
	This completes the proof.
\end{proof}


\begin{proof}[Proof of Proposition \ref{cor:dimE*}.] We may certainly assume $p(t)<\infty$.  Put $p=p(t)+\eps/3$.
	By Lemma~\ref{lem:Lambdat}, there exists $C>0$ such that if $0<\delta'<C\delta'<\delta<1/2$, then  
	$$\sum_{n=1}^\infty \sum_{z\in Z_n(\delta,\delta')} e^{-p n} e^{tS_n f_c(z)}<\infty.$$
	For each $z\in Z_n(\delta,\delta')$, let $J_{n,z}$ (resp. $J'_{n,z}$) denote the component of $T^{-n}(B(b,\delta))$ (resp. $T^{-n}(B(b,\delta'))$) containing $z$.
	As in the proof of Lemma~\ref{lem:Lambdat}, if $\delta/\delta'$ is large enough, we have that \for each $x\in J'_{n,z}$,
	$$|tS_n f_c(x)-tS_n f_c(z)|\le n\eps/3.$$
	For each $x\in E_*(\alpha; \delta, \delta')$, we can find an arbitrarily large $n$ and $z\in Z_{n}(\delta, \delta')$ such that $x\in J'_{n,z}$, and
	$$|tS_n f_c(x)-tn\alpha|\le n\eps/3,$$
	so that
	\begin{equation}\label{eqn:snzeta}
	|tS_n f_c(\zeta)-tn \alpha|\le 2n\eps/3.
	\end{equation}
	For each $N$, let $\mathcal{J}_N$ denote the collection of $J'_{n,z}$ with $n\ge N$, $z\in \bigcup_{n=N}^\infty Z_n(\delta,\delta')$ and such that (\ref{eqn:snzeta}) holds. Then
	$\mathcal{J}_N$ is a covering of $E_*(\alpha; \delta, \delta')$. Since
	$$\sum_{J'_{n,z}\in\mathcal{J}_N} e^{-n(p(t)-t\alpha+\eps)} \le \sum_{n=1}^\infty \sum_{z\in Z_n(\delta,\delta')} e^{-p n} e^{tS_n f_c(z)}<\infty,$$
	(\ref{eq:E*})  follows.
\end{proof}

\subsection{Proof of Theorem B} 

\begin{proof} (1)  For any $0<\delta'<\delta$, every point in $E(\alpha)$ eventually lands in one of the following sets:
	$$E_0(\alpha;\delta,\delta'):=E^\hyp(\alpha);\quad E_1(\alpha; \delta,\delta'):=E(\alpha)\cap\mathcal{B}(\delta,\delta');
	\quad 
	E_*(\alpha;\delta,\delta').$$
	 By Proposition \ref{lem:ptinf}, $$\log 2 \cdot \dim E_0(\alpha;\delta,\delta') \le p(t)-t\alpha$$ for each $t\in\R$.
	By Proposition~\ref{lem:dimB}, for any $\eps>0$, we can choose $\delta>\delta'>0$ so that $$
	 \log 2\cdot \dim E_1(\alpha; \delta,\delta') \le \dim (\mathcal{B}(\delta, \delta'))<\eps$$
	  and so that $\delta/\delta'$ is as large as we want. These, together with Proposition~\ref{cor:dimE*}, imply 
	$$\log 2 \cdot \dim E_c(\alpha)\le \inf_t (p(t)-t\alpha).$$
	The inverse inequality follows directly from Proposition \ref{lem:ptinf}.
	\medskip
	
	(2) Assume $\alpha(c)=-\infty$. Then $\alpha \mapsto p^*(\alpha)$ 
	is a concave  function  on $(-\infty, \beta(c))$.
	 Since $-\log 2 \le \beta(c)$ and  $p^*(-\log 2) = \log 2$ (by Theorem \ref{thm:finiteA}), we must have 
	$p^*(\alpha)  =\log 2$ for all $\alpha \le -\log 2$, so that $\dim E(\alpha) =1$.
\medskip

     (3) The set $E_c(\beta(c))$ is not empty because it  supports the Sturmian measure $\mu_{\max }$ maximizing $f_c$, which is ergodic (cf. \cite{FSS2021}, see also Introduction). To prove  $\dim E_c(\beta(c))=0$,
    it suffices to prove 
    \begin{equation}\label{eq:dimBetaC}
    \lim_{t\to+\infty} (p(t)-t\beta(c))=0.
    \end{equation} 
    Fix any $\epsilon >0$. Take $\mu_t$ such that $h_{\mu_t}+t\int f_c d\mu_t>p(t) -\epsilon$. From
    $$
    \log 2 + t \int f_c d \mu_t + \epsilon \ge p(t)\ge h_{\mu_{\max}}+t\int f_c d\mu_{\max} = t\beta(c) 
    $$
    and the fact that $\mu_{\max}$ is maximizing, we  get $\int f_c d \mu_t \to \beta(c)$ as $t \to +\infty$.
    By the upper semi-continuity of $\mu \mapsto \int f_c d\mu$, we get 
    $$
       \int f_c d \mu_{\max} \le \int f_c d \nu
    $$
    for any weak limit point $\nu$ of $\mu_t$.
    Then by
    the uniqueness of the maximizing measure, we get
    $\mu_t \to \mu_{\max}$.
    Since $\mu \mapsto h_{\mu}$ is upper semi-continuous,
    $\lim_t h_{\mu_t} \le h_{\mu_{\max}} =0$.
    Hence 
    $$
    \lim_{t\to+\infty}(p(t) - t \beta(c))
    \le \epsilon +   \lim_{t\to+\infty} t \left(\int f_c d \mu_t - \int f_c d\mu_{\max}\right)\le \epsilon,
    $$
    because the difference of two integrals are not positive. So we get
    (\ref{eq:dimBetaC}).
    
   (4)
   For $\alpha > \beta(c)$, 
   $E_c(\alpha)=\emptyset$ follows immediately  from (\ref{eq:max}). 
   Now
   suppose that $\alpha<\alpha(c)$. Then $E^\hyp(\alpha)=\emptyset$. Moreover, $p(t)-t\alpha<0$ holds when $t\to-\infty$, thus by Proposition~\ref{cor:dimE*}, 
   $E_*(\alpha,\delta,\delta')=\emptyset$ provided that $\delta/\delta'$ is large enough. Applying Proposition 
   \ref{lem:dimB}, we obtain that $\dim (E(\alpha))=0$.
\end{proof}

\end{document}